\documentstyle{amlts}
\begin{document}
\annalsline{153}{2001}
\received{February 10, 1999}
 \revised{August 25, 2000}

\startingpage{623}
\def\bye{\end{document}}
 \font\tenrm=cmr10
\input boxedeps.tex 
\SetepsfEPSFSpecial 
\HideDisplacementBoxes
\def\figin#1#2{\medbreak
$$
 {\BoxedEPSF{#1 scaled
#2}%
}%
$$
\medbreak\noindent}
\catcode`\@=11
\font\twelvemsb=msbm10 scaled 1100
\font\tenmsb=msbm10
\font\ninemsb=msbm10 scaled 800
\newfam\msbfam
\textfont\msbfam=\twelvemsb  \scriptfont\msbfam=\ninemsb
  \scriptscriptfont\msbfam=\ninemsb
\def\msb@{\hexnumber@\msbfam}
\def\Bbb{\relax\ifmmode\let\next\Bbb@\else
 \def\next{\errmessage{Use \string\Bbb\space only in math
mode}}\fi\next}
\def\Bbb@#1{{\Bbb@@{#1}}}
\def\Bbb@@#1{\fam\msbfam#1}
\catcode`\@=12

 \catcode`\@=11
\font\twelveeuf=eufm10 scaled 1100
\font\teneuf=eufm10
\font\nineeuf=eufm7 scaled 1100
\newfam\euffam
\textfont\euffam=\twelveeuf  \scriptfont\euffam=\teneuf
  \scriptscriptfont\euffam=\nineeuf
\def\euf@{\hexnumber@\euffam}
\def\frak{\relax\ifmmode\let\next\frak@\else
 \def\next{\errmessage{Use \string\frak\space only in math
mode}}\fi\next}
\def\frak@#1{{\frak@@{#1}}}
\def\frak@@#1{\fam\euffam#1}
\catcode`\@=12

\newfont{\bb}{msbm10 at 11pt}
\newfont{\bbsmall}{msbm10 at 9pt}
 \def\a{{\alpha}}
\def\S{{\Sigma}}
\def\Enod{{\partial_{\infty}E\neq S_{\infty}}}
\def\h{{\HY^3}}
\def\b{{\partial}}
\def\veh{{\overrightarrow{H}}}
\def\n{{\nabla }}
\def\O{{\Omega}}
\def\res{\vrule height .6em depth .2em width 0.05em}
\newcommand{\nor}[2]{\vert\nabla#1\vert^{#2}}
\newcommand{\R}{\mbox{\bb R}}
\newcommand{\C}{\mbox{\bb C}}
\newcommand{\HY}{\mbox{\bb H}}
\newcommand{\N}{\mbox{\bb N}}
\newcommand{\ve}[1]{\overrightarrow{#1}}


\title{The geometry of finite topology\\ Bryant surfaces}
 \shorttitle{Finite topology Bryant surfaces}   

 \twoauthors{Pascal Collin, Laurent Hauswirth,}{Harold Rosenberg}
 \institutions{Laboratoire de Math\'ematiques \'Emile Picard,\\ Universit\'e Paul Sabatier,
Toulouse, France\\
{\eightpoint {\it E-mail address\/}:  collin@picard.ups-tlse.fr}\\
\vglue6pt
Universit\'e de Marne-la-Vall\'ee, Marne-la-Vall\'ee, France\\
{\eightpoint {\it E-mail address\/}:  hauswirth@math.univ-mlv.fr}\\
\vglue6pt
Universit\'e Paris 7 Denis Diderot, Paris, France\\
{\eightpoint {\it E-mail address\/}:   rosen@math.jussieu.fr  
}}
 
\section{Introduction}  

In this paper we shall establish that properly embedded constant mean 
curvature one surfaces in $\h$ of finite topology are of finite total 
curvature and each end is regular. In particular, this implies the 
horosphere is the only simply connected such example, and the catenoid 
cousins the only annular examples of this nature. In general each annular 
end of such a surface is asymptotic to an end of a horosphere or an end of 
a catenoid cousin.

Robert Bryant  discovered a holomorphic parametrization of (simply connected) 
mean curvature one surfaces in $\h$ which can be thought of as a 
generalization of the Weierstrass representation of minimal surfaces in 
$\R^3$ \cite{Bryant}. Each (simply connected) minimal surface in $\R^3$ is 
isometric to a mean curvature one surface in $\h$ (and vice versa); R.~Bryant 
calls this the cousin of the minimal surface. This correspondence follows 
easily from Bonnet's existence theorem for surfaces in the space forms. 
This may have been R.~Bryant's motivation to seek a meromorphic Weierstrass 
type representation of mean curvature one surfaces in $\h$. 

\smallbreak {\it Definition}.
A Bryant surface is a surface in $\h$ of constant 
mean curvature~one. 
\smallbreak

The Weierstrass pair of a minimal (local) surface in $\R^3$ is a pair of 
meromorphic data $(g,\omega)$. The cousin in $\h$ has more local structure; 
in particular, one also has the hyperbolic Gauss map $G$. The surfaces are 
isometric so the metric is determined by $(g,\omega)$: 
$ds=\vert\omega\vert\left( 1+\vert g\vert^2 \right)$. However, 
the Gauss map $G$ is fundamental to the geometry of the cousin in $\h$ 
($G$ is also meromorphic on the minimal cousin in $\R^3$, but this seems 
never to have been considered).

An annular end of a finite total curvature minimal surface in $\R^3$ is 
conformally a punctured disk and the Gauss map $g$ extends meromorphically 
to the puncture. However, the Gauss map of an annular cousin in $\h$ may have 
an essential singularity at the puncture; R.~Bryant observed this for 
Enneper's minimal surface in $\R^3$ and its cousin in $\h$ \cite{Bryant}. 
An annular end in $\h$ is called regular if it is conformally a punctured 
disk and $G$ extends meromorphically to the puncture. This notion was 
introduced and developed for Bryant surfaces by M.~Umehara and K.~Yamada 
\cite{Yam-Ume2}. The idea of regular ends originated in the paper of 
R.~Schoen \cite{Sch}, where he introduced and studied regular minimal 
annular end hypersurfaces in $\R^n$.

A properly embedded minimal annular end in $\R^3$ of finite total curvature 
is asymptotic to an end of a plane or catenoid. In $\h$, a properly embedded 
Bryant annular end,  regular  and of finite total curvature, is asymptotic to an 
end of a horosphere or a catenoid cousin \cite{Tou-Ric}.

We will prove that a properly embedded Bryant annular end in 
$\h$ is of finite total curvature and regular. This is the main result of 
our work and  answers affirmatively a conjecture by M.~Umehara and 
K.~Yamada: there are no embedded irregular Bryant annular ends of finite 
total curvature \cite{Yam-Ume1}.

When the annular end is part of a properly embedded Bryant 
surface, we prove it is asymptotic to a catenoid cousin end and not a 
horosphere end (unless $M$ is equal to a horosphere). This is   
Theorem~\ref{the12}. This is quite different from properly embedded 
minimal surfaces in $\R^3$, where an annular end can be asymptotic 
to a catenoid or planar end, as in Costa's surface.

The analogous theorem for minimal annular ends in $\R^3$ is not true: 
the helicoid has an annular end of infinite total curvature. The cousin of 
the helicoid in $\h$ is not embedded. In fact, our initial motivation 
was the search for a properly embedded simply connected Bryant surface in 
$\h$, other than a horosphere (the cousin of a plane in $\R^3$). Now 
we know there is no such simply connected surface.

It is still unknown if the helicoid and plane are the only properly 
embedded minimal surfaces in $\R^3$ that are simply connected.

However, the geometry of properly embedded minimal $M$ in $\R^3$, 
of finite topology and with at least two ends, is understood: $M$ 
has finite total curvature; each annular end of $M$ is asymptotic to a 
plane or catenoid end \cite{collin}, \cite{Meeks-Rosenberg}.

There has been much important work done on the geometry of properly embedded 
annular $H$-ends in $\R^3$ and in $\h$ \cite{KKMS}, \cite{KKS}, \cite{Mee}. 
In $\R^3$, for $H\neq0$, they prove such an end is asymptotic to a 
Delauney end \cite{KKS}. Also it is proved that if $H>1$ for such an end 
in $\h$ then it is also asymptotic to a Delauney end \cite{KKMS}. In fact, 
the linking number argument of our Theorem~\ref{the9} is inspired by the 
linking number argument of \cite{KKS}. However this argument needs to be 
adapted to our situation. Essentially, because of   the noncompactness of 
horospheres, we cannot use them directly as barriers. So, we will 
construct stable surfaces with sufficiently known behavior at infinity 
and use them as comparison surfaces with horospheres.

There are examples of higher genus, mean curvature one surfaces in $\h$ of finite topology. Many such examples have been 
constructed by W.~Rossman, M.~Umehara and K.~Yamada and computer images 
indicate many of these surfaces may be embedded \cite{RUY}. In $\R^3$, 
N.~Kapouleas has constructed many properly immersed and embedded 
$H$-surfaces by desingularizing certain families of touching spheres. 
We hope  that this may be done in $\h$, by desingularizing certain 
families of touching horospheres. For example, consider the three 
horospheres intersecting in three points as in Figure~1-a. 
One should be able to attach catenoid cousin necks near 
the three singular points and show there is a Bryant surface in a  
neighborhood of this new surface by Schauder fixed point techniques. 
This surface would have genus one and three ends, each asymptotic to a 
catenoid cousin end; see Figure~1-b.

\bigbreak
\begin{center}
\BoxedEPSF{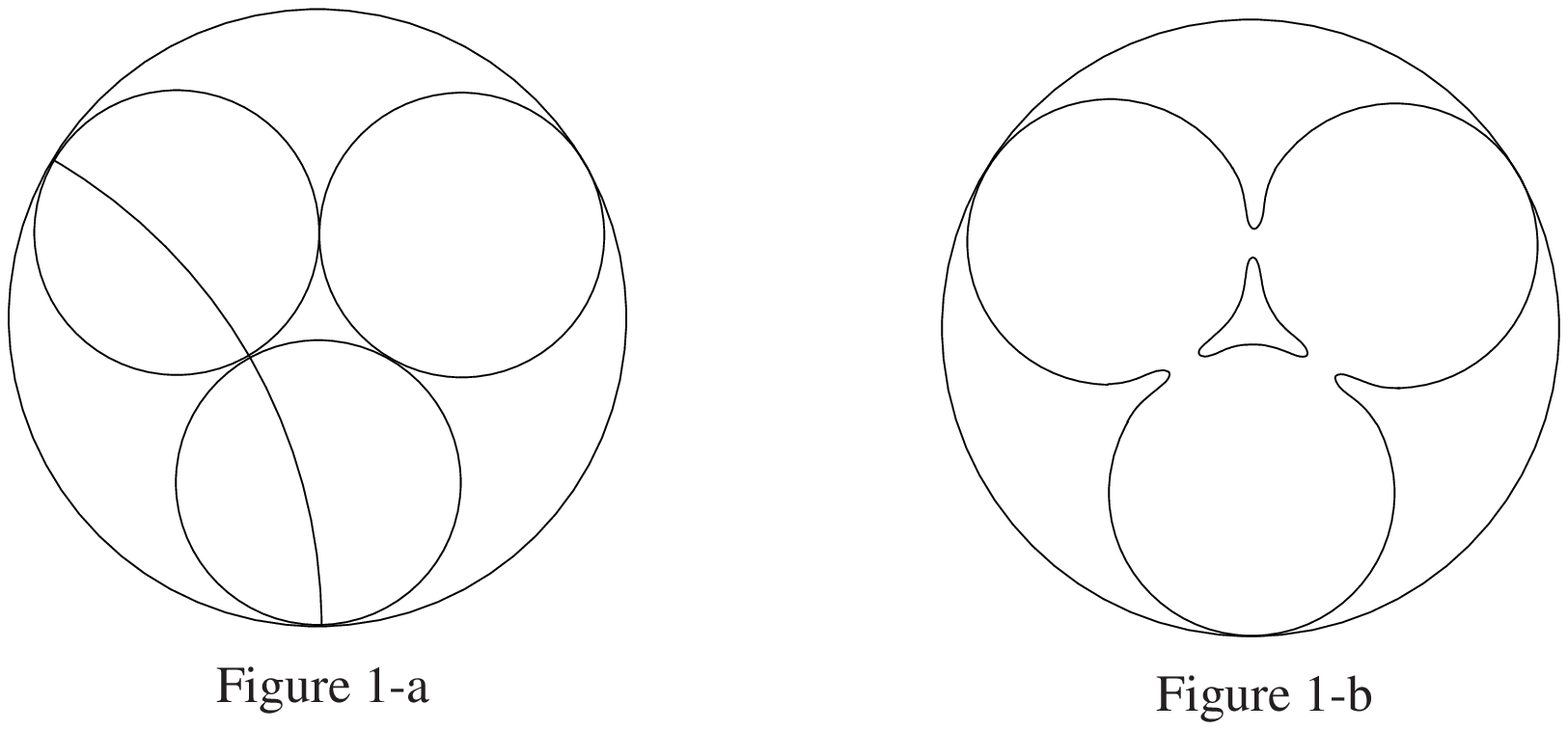 scaled 680} 
\end{center}
\smallbreak

The paper is organized as follows. In Section~2, we give a (brief) 
description of the Bryant representation; the interested reader may 
consult \cite{Bryant}, \cite{Tou-Ric} and \cite{Yam-Ume2} for a 
serious discussion.

In Section~3, we analyze the connected component of the intersection 
of a  Bryant surface in $\h$ with its tangent horosphere at a point. There 
is more structure here than the trace of a minimal surface in $\R^3$ on its 
tangent plane at a point. We describe here this trace for properly embedded 
annular ends $E$ with $E \cap H(q)$ compact, $q \in E$, and 
$\b E \cap H(q)=\emptyset$; $H(q)$ is the tangent horosphere at $q$.

In Section~4 we study properly embedded Bryant annular ends $E$, which 
are not dense at infinity. We first prove this end is regular: it is 
conformally the punctured disk and $G$ extends meromorphically to the 
puncture (Theorem~\ref{the1}). We then prove the asymptotic boundary of 
the end $E$ is precisely the limiting value of $G$ at the puncture 
(Theorem~\ref{the2}).

In Section~5 we continue the study of properly embedded annular ends 
$E$ assuming $E$ is regular. We prove $E$ then has finite total 
curvature. This is done by first proving  such an end $E$ has finite 
total curvature if it is on the mean convex side of a catenoid cousin 
(Theorem~\ref{the3}). Then we prove this end can be placed on the 
mean convex side of a catenoid cousin (Theorem~\ref{the4}). This last 
result requires an analytic theorem concerning $H=1$ graphs over 
noncompact domains (Theorem~\ref{the5}).

In Section~6 we prove the nondensity at infinity of finite topology 
properly embedded Bryant surfaces. Using the trace on horospheres, we show 
that if $M$ is dense at infinity there is a proper arc $\gamma$ on $M$ 
with $\b_\infty \gamma$ two distinct points at infinity. Using $M$ as a 
barrier we construct stable $H=1$ surfaces with boundary $\gamma$. 
Analyzing the asymptotic behavior of stable surfaces we see that such 
stable surfaces cannot exist.

\section{The Bryant representation}
\label{repres}

Let ${\cal L}^4$ be Minkowski 4-space with the Lorentzian metric 
of signature $(-,+,+,+)$. Hyperbolic 3-space can be represented as 
$$\h=\left\{ \left( t,x_1,x_2,x_3 \right) \in {\cal L}^4; \quad
\sum_{i=1}^3 x_i^2-t^2=-1,\ t>0 \right\}$$ 
with the metric induced from ${\cal L}^4$.

It is useful to identify ${\cal L}^4$ with the space of $2\times 2$ 
hermitian matrices: a point $\left( t,x_1,x_2,x_3 \right)$ corresponds to 
$$\pmatrix{t+x_3 & x_1+ix_2  \cr x_1-ix_2 & t-x_3}.$$
Notice that $\h$ is the set of such matrices of determinant one, $t>0$, 
and one has $\h=\left\{ aa^*;a \in {\rm SL}(2,\C) \right\}$, where 
$a^*=^t\overline{a}$.

Let $M$ be a simply connected Riemann surface and $F: M \rightarrow {\rm SL}(2,\C)$ a holomorphic immersion satisfying: 
$$dAdD-dBdC=0,$$
where $F=\pmatrix{A&B \cr C&D}$.

Then $f=FF^*:M \rightarrow \h$ is a conformal immersion of mean 
curvature-one. If $H \in {\rm SU}(2)$ then $f=F_1F^*_1$ where $F_1=FH$.

Conversely any mean curvature one surface in $\h$ is given locally 
by such an $F$. The reader should consult \cite{Bryant} and 
\cite{Yam-Ume2} for the details.

In the upper half-space model of $\h$, one can express the immersion 
in terms of~$F$.
\begin{eqnarray*}
\left( x_1+ix_2 \right) (z)&=&{A\overline{C}+B\overline{D} \over 
\vert C \vert^2+\vert D \vert^2}(z),\\
x_3(z)&=&{1 \over \vert C \vert^2+\vert D \vert^2}(z).
\end{eqnarray*}

The Weierstrass data of the minimal cousin in $\R^3$ are given by:
$$F^{-1}dF=\pmatrix{ g&-g^2 \cr 1&-g }\omega.$$

Then one obtains
$$g=-{B' \over A'}=-{D' \over C'},\ \omega=AC'-A'C,\hbox{ and }G={A' \over
C'}.$$

The metric induced on $M$ is  $ds=\vert \omega \vert 
\left( 1+\vert g \vert^2 \right)$.

\section{The tangent horosphere}
\label{tanhoro} 

For $M$ an immersed surface in $\h$ and $q \in M$, the tangent horosphere 
$H(q)$ of $M$ at $q$ is the horosphere tangent to $M$ at $q$ whose mean 
curvature vector at $q$ has the same direction as that of $M$ at $q$. This 
horosphere is unique when the mean curvature of $M$ at $q$ is nonzero.

Note that $H(q)$ separates $\h$ into two components. We let $H(q)^+$ denote the mean 
convex component bounded by $H(q)$ and we call it the inside of $H(q)$. 
The surfaces at a constant distance $t$ from $H(q)$ are also horospheres  
with the same point at infinity as $H(q)$ and they foliate $\h$. We denote 
this equidistant horosphere by $H_t(q)$, and for $t>0$, $H_t(q)$ will 
be inside $H(q)$ and outside $H(q)$ for $t<0$.

Now suppose $M$ is a Bryant surface properly embedded 
in $\h$. We allow $M$ to have a compact boundary since many of our results 
concern the ends of such surfaces. We also assume $M$ is not a part of a 
horosphere.

For $q \in {\rm int}(M)$, the intersection of $M$ and $H(q)$ is an analytic curve 
near $q$ with isolated singularities, and at the singularity $q$, there 
are $2k+2$ smooth branches meeting at equal angles where $k$ is an integer at 
least one. In  fact, $k$ is the same as the order of contact of the 
cousin minimal surface in $\R^3$ with its tangent plane, and $k-1$ is the 
order of $q$ as a branch point of the Gauss map $G$.

When $\b M=\emptyset$, $M$ separates $\h$ into two connected components 
since $M$ is properly embedded. We let $W$ denote the mean convex 
component bounded by $M$. When $\b M \neq \emptyset$  and is compact, we 
introduce a mean convex component $W$ as follows. It is not hard to see 
that there is an embedded compact orientable surface $\S$ such that 
$\b \S=\b M$ and $\S \cap {\rm int} (M)=\emptyset$ (take a large ball $B$ of 
$\h$, containing $\b M$, such that $M$ is transverse to $\b B$. Let 
$M_0=M\cap B$ so that $\b M_0=\b M\cup\Gamma$, where $\Gamma$ is a one-dimensional 
submanifold of $\b B$. Now $\Gamma$ bounds a compact domain $D\subset\b B$ such 
that $\overrightarrow{H}$, the mean curvature vector of $M$, points 
towards $D$ along $\Gamma$. Then $\Sigma$ can be obtained by smoothing 
$M_0\cup D$ along $\Gamma$ and displacing this slightly, in the direction 
of $\overrightarrow{H}$, keeping $\b M$ fixed). Now $\S \cup M$ separates 
$\h$ into two components and we call $W$ the component into which 
$\overrightarrow{H}$ points along $M$. We will use $W$ far from $\b M$ so that 
the choice of $\S$ is not important.

We will now derive properties of $E \cap H(q)$ where $E$ is a properly 
embedded Bryant annular end (homeomorphic to $S^1\times[1,+\infty[)$. We 
henceforth assume $\b E \cap H(q)=\emptyset$ and $E$ is topologically 
the unit disk punctured at the origin.

\specialnumber{1} \proclaim{Lemma}
\label{lemout}
Let $E_1$ be a connected component of $E-H(q)$  that is  outside 
$H(q)${\rm .} Then $E_1$ is not compact or $\b E \subset E_1${\rm .}
\endproclaim

{\it Proof}.
If this were not so, then $E_1$ would be compact and $\b E_1 \subset H(q)$. 
Consider a ``large" horosphere $H_t(q)$ outside $H(q)$ such that 
$E_1 \subset H_t (q)^+$ (so $t$ is near $-\infty$ in our notation). Then 
increase $t$: since $E_1$ is compact there will be a largest $t_0$ such 
that $H_{t_0}(q)$ touches $E_1$ at a point $x \in {\rm int}(E_1)$. But then 
$E=H_{t_0}(q)$ by the maximum principle, a contradiction; cf.\ Figure~2-a. 
\hfill\qed

\specialnumber{2} \proclaim{Lemma}
\label{lemins}
Let $E_1$ be a compact connected component of $E-H(q)$ with 
$\b E_1 \subset H(q)$ and $E_1$ inside $H(q)^+$. Let $D_1 \subset H(q)$ 
be a compact domain{\rm ,} $\b D_1=\b E_1${\rm ,} and $D_1 \cup E_1=\b Q_1${\rm ,} 
$Q_1$ a compact domain in $H(q)^+${\rm .} Then $Q_1$ is mean convex along $E_1${\rm .}
\endproclaim

{\it Proof}.
Consider a ``small" horosphere $H_t(q)$ contained in $H(q)^+-Q_1$ (so $t$ 
is near $+\infty$). When  $t$ decreases, there will be a positive $t_0$ where 
$H_{t_0}(q)$ touches $Q_1$ for the first time. The point $x$ where they touch 
is in $E_1$ and by the maximum principle, the mean curvature vector of $E_1$ 
at $x$ is the negative of that of $H_{t_0}(q)$ at $x$. So this mean curvature 
vector points into $Q_1$ and $Q_1$ is mean convex; see  Figure~2-b. 
\hfill\qed

\specialnumber{3} \proclaim{Lemma}
\label{lemcomp}
There is at most one compact component at $q$ of $E-H(q)${\rm ,} whose boundary 
is in $H(q)${\rm ;} by {\rm ``}\/at $q${\rm "} is meant a connected component of $E-H(q)$ 
containing $q$ in its closure{\rm .}
\endproclaim  

\vfill
\begin{center}
\BoxedEPSF{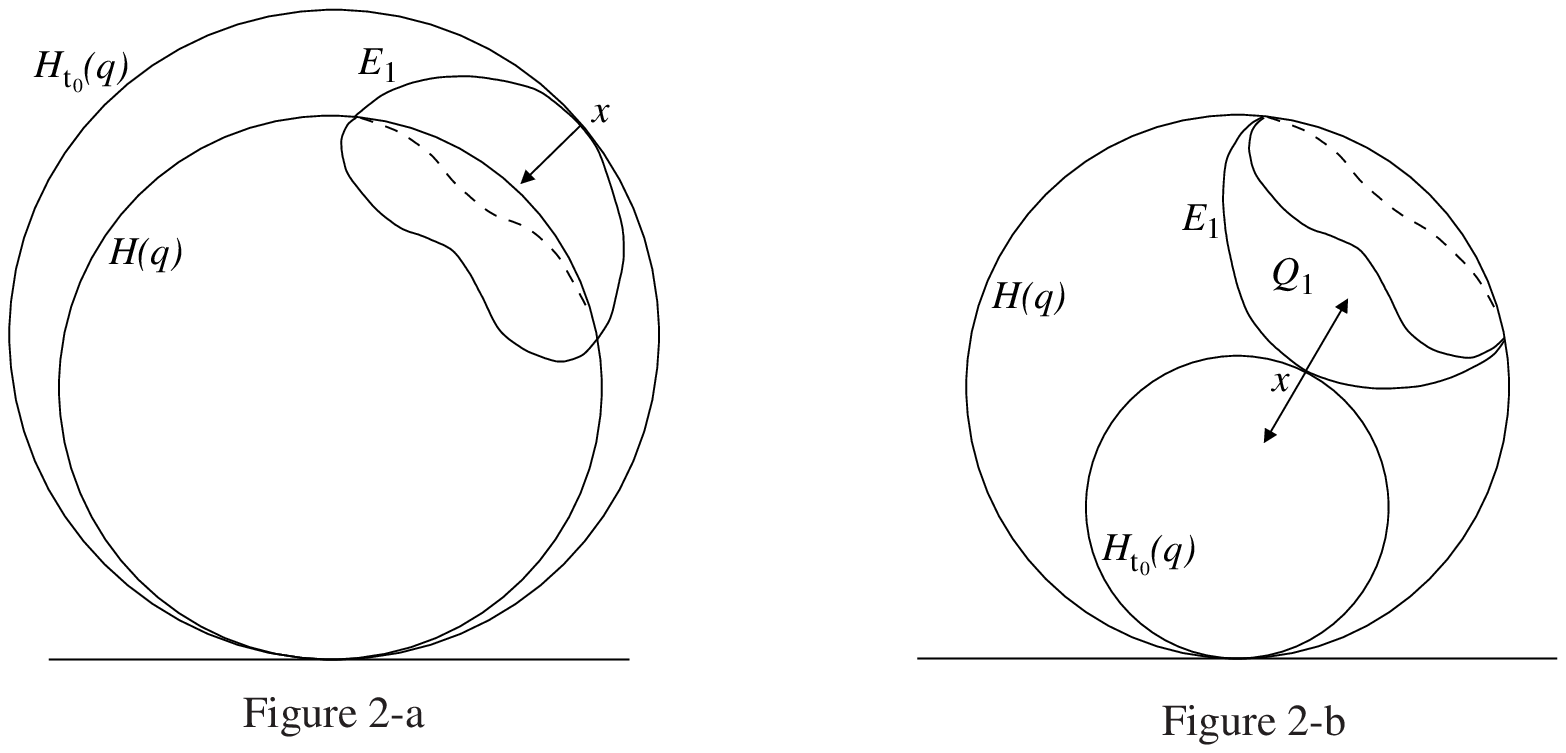 scaled 800}
\end{center}

\eject
\demo{Proof}
Suppose this fails. Then at least two components $E_1,E_2$ at $q$ are 
compact and $(\b E_1\cup \b E_2)\subset H(q)$. By Lemma~\ref{lemout}, 
we know that $E_1 \cup E_2\subset H(q)^+$. For $i=1,2$, let 
$D_i\subset H(q)$ be compact domains, $\b D_i=\b E_i$ and 
$E_i\cup D_i=\b Q_i$, where $Q_i$ is a compact domain in $H(q)^+$.

Since $E$ is a graph over $H(q)$ near $q$, the mean curvature vectors (which 
we denote by $\overrightarrow{H}$) of $E_1$ and $E_2$ point into the same 
connected component $C$ of $H(q)^+-(E_1 \cup E_2)$ near $q$. 
So $\veh(E_1)$ and $\veh(E_2)$ point into $C$ along $E_1 \cup E_2$.

If $C$ is compact, then $C$ is contained in $Q_1$ or $Q_2$,  say $Q_1$. So 
$\veh(E_2)$ points into $C \cap Q_1$. Now $E_2 \subset Q_1$ hence 
$Q_2 \subset Q_1$ so along $E_2$, $\veh(E_2)$ points to the noncompact 
component of $H(q)^+-E_2$; this contradicts Lemma~\ref{lemins}; see Figure~3.

If $C$ is noncompact, then $\veh(q)$ points into $C$, so along $E_1$, 
$\veh$ points into $C$ as well. But $\veh$ points into $Q_1$ along $E_1$ 
by Lemma~\ref{lemins}. This proves Lemma~\ref{lemcomp}. \phantom{circleof}
\enddemo 

\vglue-18pt
\figin{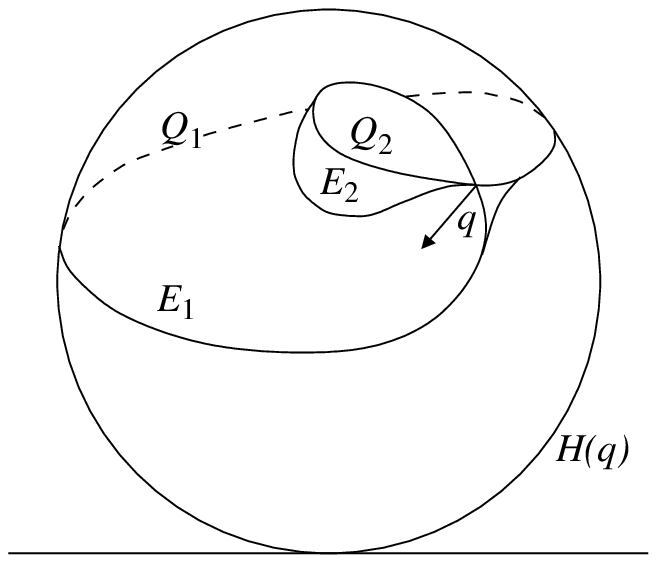}{750}
\centerline{Figure~3}

\specialnumber{4} \proclaim{Lemma}
\label{lemnum}
The number of connected components of $E-H(q)$ at $q$ is at least three{\rm .} If 
equality holds then the order of contact $k$ of $E$ with $H(q)$ is one{\rm .}
\endproclaim

\vglue-9pt
{\it Proof}.
Let $2k+2$ denote the number of branches of $E\cap H(q)$ at $q$. Let 
$f:D\rightarrow E$ be a local parametrization of a neighborhood of $q$ 
on $E$ by a disk of $\C$ so that $f(0)=q$ and lines passing through $0$ 
of slope an integral multiple of $\pi \over k+1$ are sent to $E\cap H(q)$. 
Let $A_i$ be the sector of $D$ defined by 
$\{{(i-1)\pi \over k+1}<\arg(z)<{i\pi \over k+1}\}$ and $B_i=f(A_i)$.
We know that the $B_i$ are alternatively in $H(q)^+$ and $H(q)^-$ 
around $q$.

If $B_1$ and $B_3$ are not in the same component of $E-H(q)$, then by the 
observation above, $B_2$ yields a third component so that Lemma~\ref{lemnum} 
is true.

If $B_1$ and $B_3$ are in the same component of $E-H(q)$, we can construct 
a cycle $\alpha_{13}$ on $E$ as follows: let $a_1$, $a_3$ be two points of 
$A_1$, $A_3$ respectively, and $\beta$ a path in $E-H(q)$ from $f(a_1)$ to 
$f(a_3)$; $\alpha_{13}=\beta\cup\beta'$ where $\beta'$ is the image by $f$ 
of the line segment from $a_1$ to $0$, followed by the line segment from 
$0$ to $a_3$. Now $\alpha_{13}$ meets $H(q)$ exactly at~$q$.

If $B_2$ and $B_4$ were in the same component, we could find a cycle 
$\alpha_{24}$ on $E$ which meets $\alpha_{13}$  in a single point, which 
is impossible since the genus of $E$ is zero. Thus we get at least three 
components in this case as well.

Now we study the case of equality. In this case there is only one component 
at $q$ in either $H(q)^+$ or $H(q)^-$. Then we can assume that all the 
$B_i$ for $i$ odd are in the same global component of $E-H(q)$. If 
$k\geq 2$ this means that we can construct cycles $\alpha_{15}$ (with 
$A_1$ and $A_5$), and $\alpha_{35}$ (with $A_3$ and $A_5$) exactly as we 
constructed $\alpha_{13}$. As before, these three cycles separate the 
components of $B_2$, $B_4$ and $B_6$ in $E-H(q)$. Hence we obtain at least 
four components at $q$ in this case. This completes the proof of 
Lemma~\ref{lemnum}.
\hfill\qed\medbreak 

We define $\S(q)$ to be the connected component of $q$ in $E\cap H(q)$ 
and we assume $\S(q)$ is compact in the rest of this section.

\specialnumber{5} \proclaim{Lemma}
\label{lemtroi}
$E-H(q)$ has exactly three components at $q${\rm :} $B${\rm ,} $E_1$ and $E_2${\rm .} The first{\rm ,}  $B${\rm ,} 
is compact and contains $\b E${\rm ;} $E_1$ is compact with boundary in $H(q)${\rm ;} 
$E_2$ is noncompact {\rm (}\/recall the assumption that $\b E\cap H(q)=\emptyset$ 
and $\S(q)$ is compact\/{\rm ).}
\endproclaim

{\it Proof}.
First we prove $E-H(q)$ has at most one noncompact component at $q$. 
This is immediate if $E\cap H(q)$ is compact (in fact $E\cap H(q)$ will 
always be compact until we arrive at Theorem~\ref{the7} of this paper. 
There we will need to work with the weaker hypothesis: $\S(q)$ is 
compact). Let $F$ be the noncompact component of $E-\S(q)$. For 
$\varepsilon>0$, $\varepsilon$ small, the points of $F$ a distance 
$\varepsilon$ from $\S(q)$ form a compact curve $C$ disjoint from 
$H(q)$. Since $E$ is an annulus, $C$ is in fact connected. This Jordan 
curve $C$ separates $\S(q)$ from the puncture.

With the notation  of Lemma~\ref{lemnum}, if two $B_i$ and $B_j$ 
($i\neq j$) are in noncompact components, then a curve in $E-H(q)$ 
from $q$ to the puncture starting in $B_i$ (and $B_j$) meets $C$. So 
$B_i$ and $B_j$ are in the same global component.

Then by Lemma~\ref{lemnum}, there are at least two components at $q$ 
of $E-H(q)$ that are compact. At most one may contain $\b E$, and 
the others are compact with boundary in $H(q)$ (and there is at least 
one). But such a component is in $H(q)^+$ by Lemma~\ref{lemout} and 
is unique by Lemma~\ref{lemcomp}. We call $E_1$ this unique component. 
Then the other compact component necessarily contains $\b E$ -- we 
call it $B$ --, and the third component ($E_2$, say) is noncompact.

Finally $E-H(q)$ has exactly three components whose closure contains $q$ 
and by Lemma~\ref{lemnum}, the order of contact $k$ of $E$ with $H(q)$ 
is one.
\hfill\qed

\specialnumber{6} \proclaim{Lemma}
\label{lem6}
There are two possibilities\/{\rm :} 

 {\rm --} $\b E_1=C_1$ is a Jordan curve on $H(q)${\rm ,}

\noindent or

 {\rm --} $\b E_1=C_1\cup C_2$ is a figure eight\/{\rm ;} the union of two 
Jordan curves $C_1$, $C_2$ on $H(q)$ meeting at $q${\rm .}
\endproclaim

\demo{Proof}
First notice that $\b E_1$ contains no cycle $c$ disjoint from $q$. To 
see this, let $\gamma$ be a path in $E_2\cup B\cup \{q\}$ going from 
$\b E$ to the puncture and meeting $H(q)$ exactly at $q$. The cycle 
$c$ does not meet $\gamma$ so $c$ bounds a compact domain $D$ in $E$, 
$D\cap\b E=\emptyset$ and $D\neq E_1$ (since $q\notin c=\b D$). 
But $D$ would contain a compact domain outside $H(q)$ with boundary 
in $H(q)$ contradicting Lemma~\ref{lemout}. This proves each cycle 
in $\b E_1$ meets $q$, $E_1$ is a disk and $\b E_1-\{q\}$ is an 
embedded curve.

We know that we have locally at $q$ exactly four components $B_1$, $B_2$, 
$B_3$, $B_4$; $B_1$ and $B_3$ in $H(q)^+$. Assume $B_1$ is in $E_1$ and 
$B_3$ is not in $E_1$. Then $q$ is not a double point of $\b E_1$ and 
$\b E_1=C_1$ is  a Jordan curve on $H(q)$. On the contrary if $B_3$ is also 
in $E_1$, then $\b E_1=C_1\cup C_2$; $C_1$ and $C_2$ Jordan curves 
meeting exactly at $q$; i.e., $\b E_1$ is a figure eight.
\enddemo

\vglue6pt

Let $D \subset H(q)$, $Q_1 \subset H(q)^+$ be such that 
$\b Q_1=E_1 \cup D$, where  $Q_1$ is compact.

\vglue12pt

\specialnumber{1} \proclaim{Proposition}
\label{protang}
With the   notation of Lemma~{\rm \ref{lemtroi}}\/{\rm :}

\smallbreak
 
If $\b E_1=C_1${\rm ,} then $\b E \subset Q_1$ and every divergent path starting 
at $x \in \b E$ must intersect $E_1 \cup H(q)$ at a point other than $x${\rm .}

 \smallbreak
If $\b E_1=C_1 \cup C_2${\rm ,} then every path starting at $\b E${\rm ,} staying in 
$W${\rm ,} and diverging in $\h$ must intersect $H(q)${\rm ;} $Q_1$ separates $W${\rm .} 
Moreover $\b E$ is outside $H(q)${\rm ,} and $C_1$ and $C_2$ are each 
homologous to $\partial E$ on $E${\rm .}
\endproclaim 


{\it Proof}.
The two possibilities are given by Lemma~\ref{lem6}. If $\b E_1=C_1$, then 
$D$ is the disk of $H(q)$ bounded by $C_1$ and $D \cup E_1=\b Q_1$. At 
$q$, $\veh(q)$ must point into $Q_1$ (since it does so at points of $E_1$ 
near $q$) by Lemma~\ref{lemins} so that $E_2$ or $B$ is inside $Q_1$ near $q$. 
It cannot be $E_2$ since $E_2$ is noncompact (the puncture is in $E_2$) 
and $E_2$ is properly embedded. Thus $B$ is inside $Q_1$ near $q$. Thus, 
$B \subset Q_1$ and, in particular, $\b E \subset Q_1$. Now any path 
starting at a point of $\b E$ and diverging in $\h$, must intersect 
$\b Q_1=E_1 \cup D$; see Figure~4. This proves the first 
assertion of the proposition.

Now suppose $\b E_1=C_1 \cup C_2$. Again $\veh(q)$ points into $Q_1$ and 
$Q_1$ is mean convex along $E_1$. By Lemma~\ref{lemout}, it is clear 
that $C_1$  and $C_2$  are not homologous to zero in $E$, hence each 
$C_1$, $C_2$ is homologous to $\partial E$ in $E$. Let $D_1$ be the disk 
of $H(q)$ bounded by $C_1$, $D_2$ bounded by $C_2$. We have 
${\rm int}(D_1)\cap {\rm int}(D_2)=\emptyset$ or one disk is contained in the other. 
This latter case is impossible. For if $D_2\subset D_1$, we have 
$\b Q_1\subset D_1$. Locally at $q$, $E_1$ is a graph over two opposite 
sectors of $H(q)$ (the projection of the $B_1$ and $B_3$ of 
Lemma~\ref{lem6}). By the local structure of $E_1$ near $q$, the two 
complementary sectors are in $\b Q_1$, hence in $D_1$. As $D_1$ is 
a disk on $H(q)$, at least one projection of $B_1$ or $B_3$, $B_3$ say, 
must be in $D_1$. Now $E_1$ is a disk with two points on the boundary  \pagebreak
identified at $q$. Then $\b E_1-\{q\}=C_1\cup C_2$ hence one of the two 
boundary arcs of $B_3-\{q\}$ (near $q$) is in $C_1$ and the other in $C_2$. 
But then, points of $C_1$ are in the interior of $D_1$, which is a 
contradiction.

Then $\b Q_1=E_1 \cup D_1 \cup D_2$ so $Q_1$ separates $\b E$ from infinity 
in $W$: any path starting at $\b E$, in $W$ and diverging in $\h$, must pass 
through $Q_1$. Moreover $B$ and $E_2$ are outside $H(q)$, in particular 
$\b E$ is outside $H(q)$; see Figure~5.
\hfill\qed
\demo{{R}emark {\rm 1}}
Assume $\S(q)$ is compact and $E$ is transverse to $\S(q)-\{q\}$. Then 
$\S(q)$ is a figure eight, the union of two Jordan curves $C_1$, $C_2$ 
meeting at~$q$.
\enddemo
\vglue-12pt
\figin{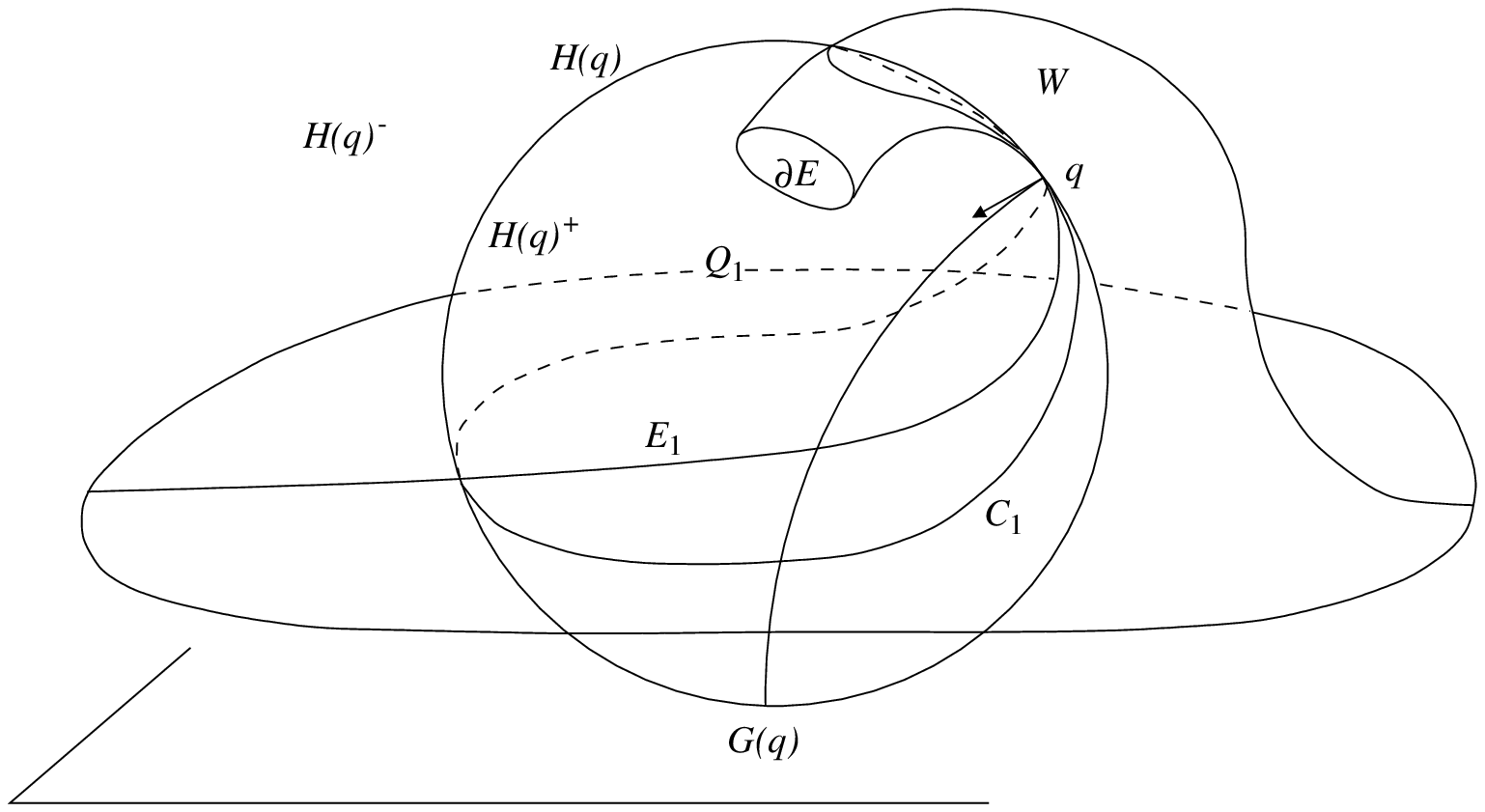}{700}
\centerline{Figure~4} 
 \vfill
\figin{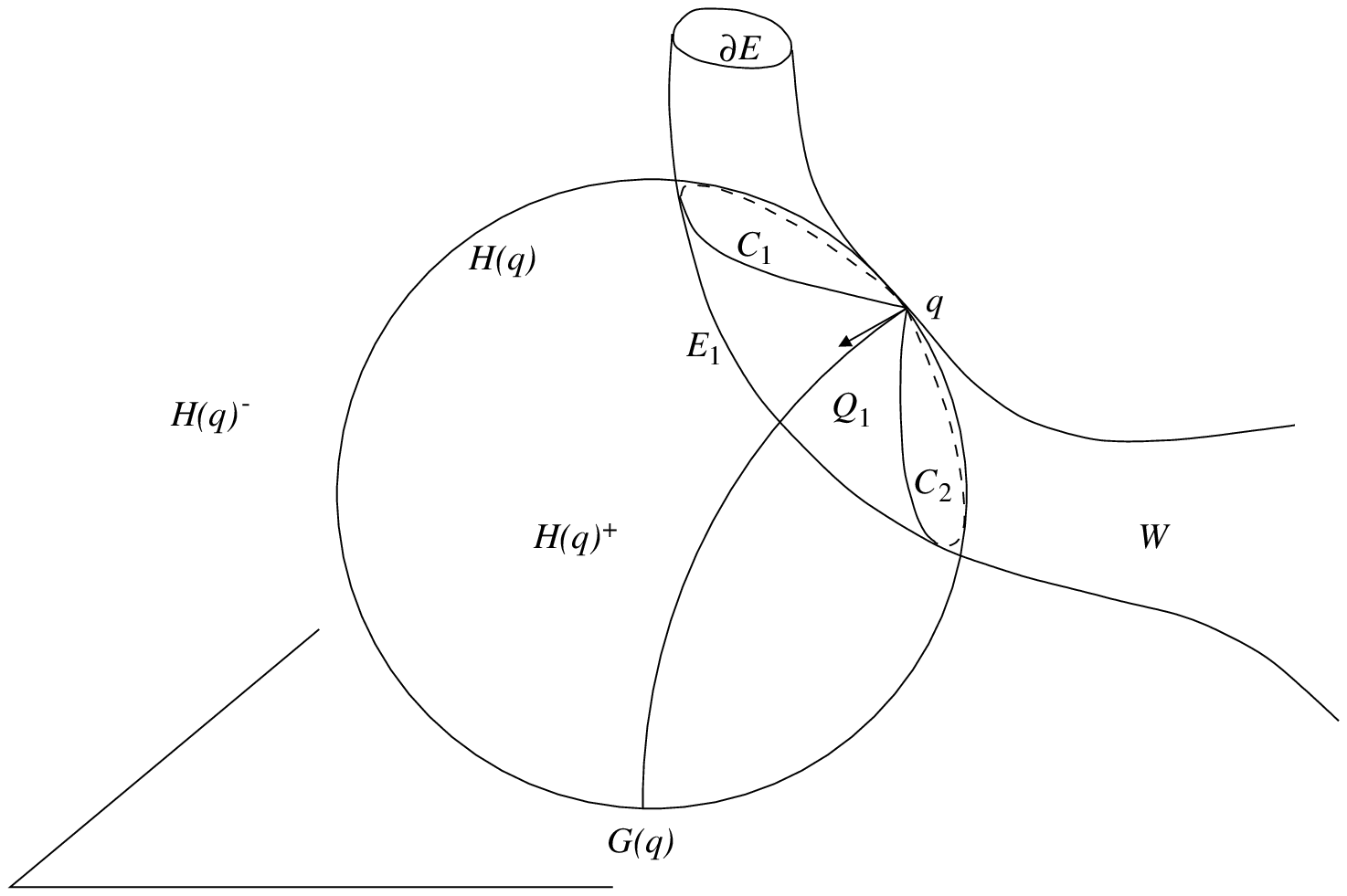}{700} 
\centerline{Figure~5} 
  \eject

By the transversality hypothesis and the local structure at $q$, $\S(q)$ 
consists of two analytic curves meeting at equal angles at $q$. It is then 
a figure eight. Note that in the first case of Proposition~\ref{protang}, 
$C_2$ together with $\b E$ bounds $B$ (other curves in $\b B$ would 
give rise to a compact component outside H(q)).

\specialnumber{1} \proclaim{{C}orollary}
\label{cortang}
Let $E$ be a properly embedded Bryant annular end{\rm ,} $q \in E$ with 
$E \cap H(q)$ compact and disjoint from $\b E${\rm .} Then if 
$\b E \subset H(q)^+$ every divergent path starting at $x \in \b E$ 
must intersect $E \cup H(q)$ at a point other than $x${\rm ;} if 
$\b E \subset H(q)^-$ then every divergent path starting at 
$x \in \b E$ and staying in $W${\rm ,} must intersect $H(q)${\rm .}
\endproclaim

\section{The regularity and asymptotic boundary\\ of an annular end not %
dense at infinity}

\specialnumber{1} \proclaim{Theorem}
\label{the1}
Let $E$ be a properly embedded Bryant annular end{\rm .} If 
$\b_{\infty}E \neq S_{\infty}${\rm ,} then $E$ is conformally a punctured 
disk and the hyperbolic Gauss map $G$ extends meromorphically to the 
puncture {\rm (}\/i.e.{\rm ,} $E$ is regular\/{\rm ).}
\endproclaim

\demo{Proof}
We will now work in the upper half-space model of $\h$ with 
$S_{\infty}=\{x_3=0\}\cup\{\infty\}$. Since the asymptotic boundary of 
$E$ is closed and not $S_{\infty}$, we can assume 
$E\subset B_R=\left\{x_1^2+x_2^2+x_3^2<R^2,x_3>0\right\}$.

First we will show that $G$ is bounded on some subend of $E$. If not, then 
for some $q_n\in E$, diverging on $E$, we would have $\vert G(q_n) \vert 
\rightarrow \infty$. Since $\b E$ is compact, we have $x_3\res_{\b E} 
\geq \delta>0$ for some $\delta$. Choose $n$ sufficiently large so that 
$H(q_n) \cap B_R$ is below $x_3=\delta$. This is possible since $x_3(q_n) 
\rightarrow 0$; cf.\ Figure~6. However, $\b E \subset H(q_n)^+$ 
and we can find a path from $\b E$ to $G(q_n)$ which does not intersect 
$E \cup H(q_n)$ except at its endpoint (choose a path from a point of 
$\b E$ to a point of $\b B_R$, not meeting $E$; then choose $n$ big enough 
so that $H(q_n)$ is below this path, and then continue to $G(q_n)$).

This contradicts Corollary~\ref{cortang}, and so $G$ is bounded.

\medbreak
 \figin{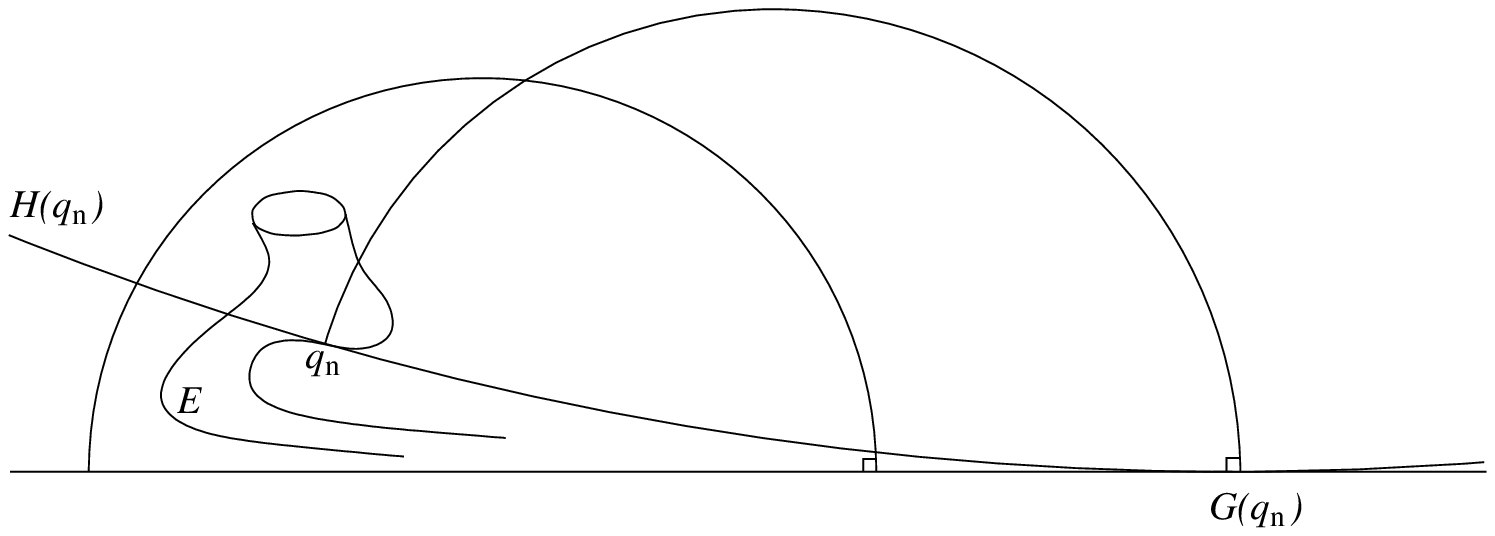}{650}
 \centerline{Figure~6}
 \medbreak

To prove Theorem~\ref{the1}, it suffices to prove that $E$ is conformally 
the punctured disk.

We will prove this by constructing a complete metric on $E$ of the form 
$d \sigma=\lambda \vert dz \vert$ where $\lambda$ is the module of a 
holomorphic function on $E$ (R.~Osserman \cite{Oss}).

Let $\widetilde{E}$ be the universal cover of $E$, so that 
$F:\widetilde{E} \rightarrow {\rm SL}(2,\C)$ is holomorphic, 
$F^{-1}dF=\pmatrix{g&-g^2 \cr 1&-g \cr}\omega$, where $(g,\omega)$ are 
the Weierstrass data and $\Psi:\widetilde{E} \rightarrow \h$, 
$\Psi=F^t\overline{F}$, defines the immersion of $E$ in $\h$. The metric 
$ds=\vert \omega \vert \left( 1+\vert g \vert^2 \right)$, and 
the meromorphic map $G$ are well-defined on $E$.

There is a dual immersion (with $H=1$) $F^{\#}:\widetilde{E} \rightarrow {\rm SL}(2,\C)$ defined by $F^{-1}:\widetilde{E} \rightarrow {\rm SL}(2,\C)$, 
introduced by M.~Umehara and K.~Yamada \cite{Yam-Ume1}.

The Weierstrass data $\left(g^{\#},\omega^{\#}\right)=
\left(G,-{g' \over G'}\omega \right)$, and $\Psi^{\#}:\widetilde{E} 
\rightarrow \h$ is $(F^{-1})^t(\overline{F^{-1}})$. This immersion need 
not define an immersion of $E$ in $\h$ but the metric $ds^{\#}$ is well-defined and nonsingular since $\Psi^\#$ is an immersion:
$$ds^{\#}={\left( 1+\vert G \vert^2 \right) \over \vert G' \vert} 
\vert g' \vert \vert \omega \vert.$$
In particular, $g' \omega / G'$ is a nonvanishing holomorphic form.

Since $\vert G \vert$ is bounded, the metric
$$d\sigma={\vert g' \vert \over \vert G' \vert} \vert \omega \vert$$
will be complete if $ds^{\#}$ is complete. Thus it suffices to prove 
$ds^{\#}$ is complete on $E$.

Let $\gamma$ be a divergent path on $E$, which is proper so that $\gamma$ 
diverges in $\h$. Now in the Lorentzian model of 
$$\h=\left\{ \left(x_1,x_2,x_3,t\right) \in 
{\cal{L}}^4;x_1^2+x_2^2+x_3^2-t^2=-1,t>0 \right\},$$ 
the path $\gamma$ diverges so that $t(\gamma) \rightarrow \infty$.

Writing $F:\widetilde{E} \rightarrow {\rm SL}(2,\C), F=\pmatrix{A&B \cr C&D }$, 
we have $2t=\vert A \vert^2+\vert B \vert^2+
\vert C \vert^2+\vert D \vert^2$. The dual immersion 
$F^{\#}=F^{-1}=\pmatrix{ D&-B \cr -C&A}$ so that $t^{\#}=t$.

In particular $t^{\#}(\gamma) \rightarrow \infty$ as well, and $\gamma$ 
diverges in $\h$ on the dual surface. Since $ds^{\#}$ is the induced 
metric on the dual surface from its immersion in $\h$, the $ds^{\#}$ 
length of $\gamma$ is infinite. This proves Theorem~\ref{the1}.
\enddemo

\demo{{R}emark {\rm 2}}
The metric $ds^{\#}$ gives information on values of the Gauss map $G$. 
Zu-Huan Yu has proved $G$ is constant if $G$ misses more than four points 
\cite{Yu}; he proves more generally that $ds^{\#}$ is complete. We have 
proved $G$ can miss at most three points when $M$ has finite total 
curvature (and $M$ is not a horosphere) \cite{CHR}.
\enddemo 

\specialnumber{2} \proclaim{Theorem}
\label{the2}
Let $E$ be a properly embedded Bryant annular end{\rm .} If $E$ is conformally 
the punctured disk $D^*$ and $G$ extends meromorphically to the puncture{\rm ,} 
then $\b_{\infty}E=G(0)$ {\rm (}\/the value of $G$ at the puncture\/{\rm ).}
\endproclaim

\demo{Proof}
As in the proof of Theorem~\ref{the1}, we work in the upper half-space 
model and assume $G(0)$ is the point at infinity. First observe that 
$G(0)\in\b_{\infty}E$. For otherwise -- since $E$ is proper and 
$\b_{\infty}E$ is closed on $S_{\infty}$ --, $E$ would be contained 
in some half space (a complement of a neighborhood of $\infty$) 
$B_R=\{x_1^2+x_2^2+x_3^2<R^2,x_3>0\}$. Then (as in the proof of 
Theorem~\ref{the1}) $G$ must be bounded; a contradiction.

Now suppose $E$ accumulates at another point at infinity which we may 
assume $\sigma=(0,0,0)$. Let $C_T$ be the cylinder 
$\{x_1^2+x_2^2 \leq T^2,0<x_3<T\}$. There are points of $E$ in $C_T$ 
for all $T>0$. As $q$ diverges on $E$, towards $\sigma$, $G(q)$ tends 
to infinity. The geodesic normal to $E$ at $q$ is a half circle meeting 
$x_3=0$ at two points, one point close to $\sigma$ (close in the metric 
$dx_1^2+dx_2^2$) and the other point $G(q)$ that is ``far" from $\sigma$. 
Thus the mean curvature vector $\veh(q)$ of $E$ at $q$, tends to a 
vertical vector pointing up. As $x_3(q)$ decreases this vector becomes 
more vertical.

Now choose $T$ sufficiently small that $q \in C_T$ implies the angle 
between $\veh(q)$ and $\overrightarrow{e_3}=(0,0,1)$ is less than 
$\pi/8$.

Then for $q \in E \cap C_T$, the vertical segment going down from $q$ to 
$S_{\infty}$ does not meet $E$ again, since $E$ bounds a mean convex 
domain $W$ (this makes sense since $C_T$ can be  chosen far from $\b E$). 
So $E \cap C_T$ is a vertical graph $u$ over a (possibly disconnected) 
planar domain.

Now we prove that for $T$ sufficiently small, $E \cap C_T$ is a vertical 
graph over the whole base of $C_T$: $x_1^2+x_2^2<T^2$. Since there can 
be no points of $E$ below this graph, this contradicts 
$\sigma \in \b_{\infty}E$.

We now make useful gradient estimates for this graph $u$ at 
$q \in E \cap C_T$ in the Euclidean metric.

Consider the vertical plane $Q$ containing the unit normal vector 
$\overrightarrow{n}$ to $E$ at $q$. $G(q)$ is also in this plane and we 
have   Figure~7 in the plane $Q$.

\begin{center} \BoxedEPSF{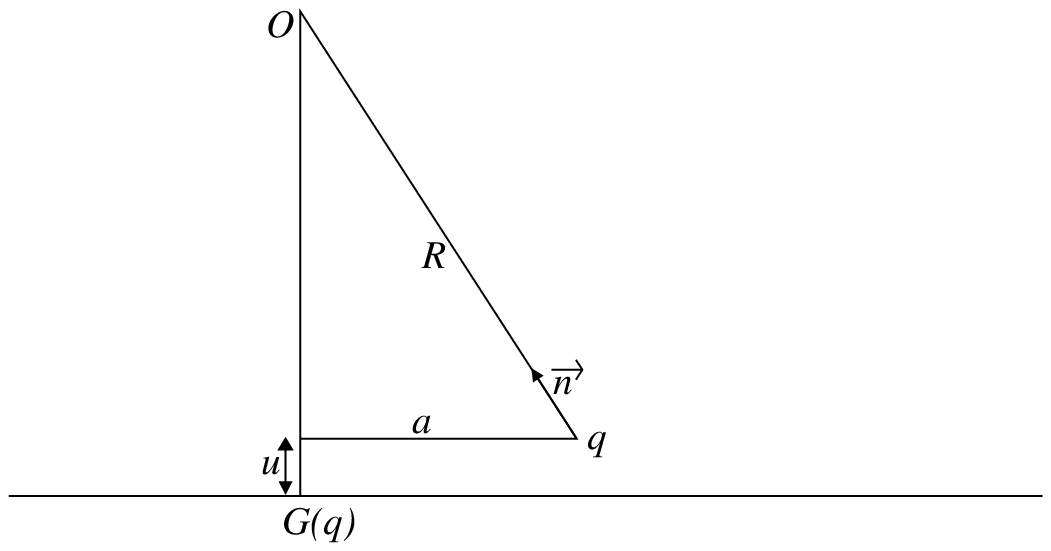 scaled 700}
\end{center}

\centerline{Figure~7} 
\medbreak

Here $O$ is the center of $H(q)$ and $R$ is the radius of $H(q)$. Then 
$$\overrightarrow{n}={1 \over W}(-u_{x_1},-u_{x_2},1), \quad
W=\sqrt{1+\vert \nabla u \vert^2}.$$

We have $\displaystyle{a=\left| {R \over W}(-u_{x_1},-u_{x_2}) \right|=
{\vert \nabla u \vert \over W}R}$, and $a^2=R^2-(R-u)^2=\break u(2R-u)$. Hence 
\begin{eqnarray*}
  {a^2 \over R^2}={\vert \nabla u \vert^2 \over W^2}&=&
 {u(2R-u) \over R^2} ,\\
 {\vert \nabla u \vert^2 \over uW^2}&=&{2R-u \over R^2}.
\end{eqnarray*}
Thus the horizontal component of $\overrightarrow{n}$ has length 
$\displaystyle l={2au \over a^2+u^2}$, and the vertical component length 
$\displaystyle t={a^2-u^2 \over a^2+u^2}$.

Now $G(0)=\infty$ and so for any large $b>0$ we can assure that $a>b$ in $C_T$ 
for $T$ small enough and 
$$\vert \nabla u \vert={l \over t}={2\left({u\over a}\right)\over 1-
\left({u\over a}\right)^2}\leq4\left({u\over b}\right).$$
Then the auxiliary function $v=\ln u$ has bounded gradient.

Starting with $x_3(q)$ small with respect to $T$, we have $v\leq\ln(T/2)$ 
on the base of $C_T$, and $E$ does not leave $C_T$ at the top $\{x_3=T\}$. 
Moreover $v$ is never $-\infty$ hence $E$ never reaches $\{x_3=0\}$ in 
$C_T$. Thus $E \cap C_T$ is a graph over the base of $C_T$ and 
Theorem~\ref{the2} is proved.
\enddemo

Theorems~\ref{the1} and \ref{the2} immediately imply:

\specialnumber{2} \proclaim{{C}orollary}
\label{cororeg}
Let $E$ be a properly embedded Bryant annular end{\rm .} If 
$\b_{\infty}E \neq S_{\infty}${\rm ,} then $E$ is regular and 
$\b_{\infty}E$ is the limiting value of $G$ on $E${\rm .}
\endproclaim

\section{Finite total curvature of nondense annular ends}

\specialnumber{3} \proclaim{Theorem}
\label{the3}
Let $E$  be a properly embedded Bryant annular end{\rm .}  If $E$ is on the mean 
convex side of a catenoid cousin end{\rm ,} then $E$ has finite total curvature{\rm .}
\endproclaim

{\it Proof}.
First we make precise ``the mean convex side." The ends of the family of 
catenoid cousins can be written as graphs (in the upper half-space model) 
over domains at infinity: $x_1^2+x_2^2 \geq r_0^2$, $x_3=0$. These ends are 
asymptotically ${1 \over r^{\alpha}}$, $\alpha>-1$ and $r^2=x_1^2+x_2^2$. 
Let $C_\alpha$ be such a catenoid cousin end and extend $C_\alpha$ to an 
embedded surface with no boundary by attaching the horizontal disk along 
$\b C_\alpha$. The mean convex side of $C_\alpha$ is then the component to 
which $\veh$ points along $C_\alpha$ (here $\veh$ is pointing up). So our 
hypothesis on $E$ is that $E$ is contained in this mean convex side of 
$C_\alpha$. Clearly the catenoidal ends are ordered by $\alpha$ and we 
can assume $\alpha>0$.

Since $\b_{\infty}E$ is the point at infinity, Theorem~\ref{the1} 
applies and we know $E$ is conformally a punctured disk 
$D^*=\{0<\vert z \vert\leq 1\}$, and $G$ extends meromorphically to $0$. 
Parametrize so that $G(z)={1 \over z^p}$ for some integer $p \geq 1$.

\smallbreak
The end $E$ is determined by $F:\widetilde{E} \rightarrow {\rm SL}(2,\C)$, 
$F=\pmatrix{A&B \cr C&D}$, with $C=z^{\nu}f$, $f$ holomorphic in $D^*$, 
and similar representations for $A$, $B$ and $D$ (this is proved in 
Lemma~\ref{periode}, following the present proof).

We know that $x_3={1 \over \vert C \vert^2+\vert D \vert^2}$. 
Suppose $C=z^{\nu}f$ and $f$ has an essential singularity at $0$. 
Then for some sequence $z_n \rightarrow 0$, we have 
$$\vert C(z_n) \vert^2 \geq {1 \over \vert z_n \vert^{(p+1)\alpha}}.$$

Let $q_n$ be the point on $E$ corresponding to $z_n$. Since $E$ is above 
the catenoid cousin $C_\alpha$: 
$$r(q_n)^{\alpha} \geq {1 \over x_3(q_n)},$$
so by the previous inequality for $C(z_n)$, we conclude 
$r(q_n) \geq {1 \over \vert z_n \vert^{(p+1)}}$, and for any integer 
$k>1$, and $n$ sufficiently large: 
$$r(q_n) > {k \over \vert z_n \vert^p}.$$

That is, the horizontal (Euclidean) distance from the point $q_n$ to the 
point $s=(0,0,x_3(q_n))$ is at least ${2 \over \vert z_n \vert^p}$. 
Observe that $d(q_n,G(q_n))$ is at least $d(q_n,s)-d(G(q_n),s)$ where 
$d$ denotes the horizontal Euclidean distance.

Let $l$ be the horizontal disk of diameter ${2 \over \vert z_n \vert^p}$, 
centered at the point $p=(G(q_n),$ $x_3(q_n))$. Since the horizontal 
distance from $G(q_n)$ to $(0,0)$ is ${1 \over \vert z_n \vert^p}$, 
the disk $l$ is in the interior of $H(q_n)^+$; see Figure~8.

Now the origin is under one of the boundary points of $l$. Observe that 
the catenoid cousin $C_\alpha$ is above the segment $[p,s]$ on $l$, since 
the height of $C_\alpha$ at $G(q_n)$ is asymptotically 
${1 \over \vert G(q_n) \vert^{\alpha}}=\vert z_n \vert^{p\alpha}$ and 
$x_3(q_n) \leq \vert z_n \vert^{(p+1)\alpha}$. Since the graph of 
$C_\alpha$ is monotone decreasing with $r$, the segment $[p,s]$ is below 
$C_\alpha$. Also, $E$ is above $C_\alpha$ so that $[p,s]$ is disjoint from $E$.

Moreover let $N$ be a compact embedded surface with boundary the boundary 
of $E$ so that $N \cup E$ is an embedded surface. $N$ can be chosen above  the union of 
the catenoid cousin $C_\alpha$ and the flat disk capping off $C_\alpha$. 
Then exactly as in Section~\ref{tanhoro}, $N \cup E$ separates the ambient 
space so one can find a path $\gamma$ from $s$ to $\b E$ which meets the 
$N \cup E$ only at the endpoint (first vertical, then a fixed path). The $k$ of the above inequality can be chosen large enough
so that this  path, together with the boundary of $E$, is inside $H(q_n)$.

\phantom{space}
\begin{center} \BoxedEPSF{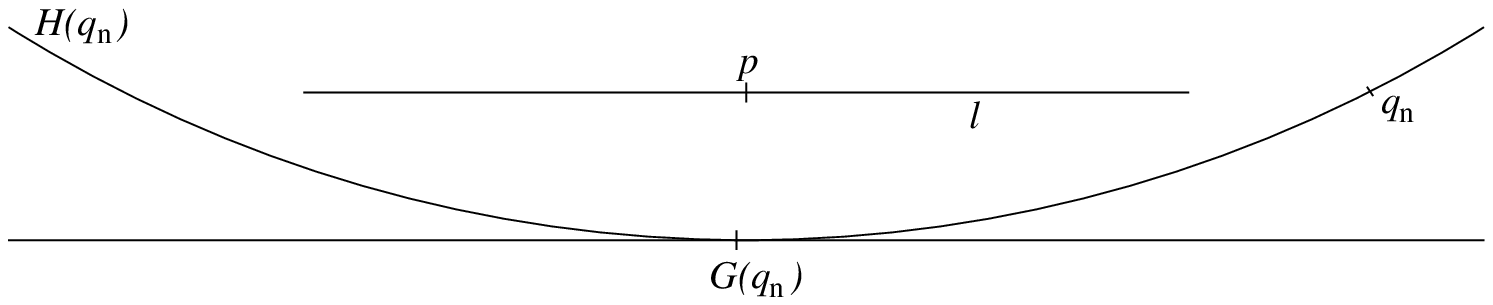 scaled 650}
\end{center}
 \centerline{Figure~8}
 \medbreak

But $\gamma$ together with [p,s] can be extended to a divergent path 
disjoint from $E \cup H(q_n)$, by going down vertically to $G(q_n)$ from 
$p$. This divergent path from $\b E \subset H(q_n)^+$ does not meet 
$H(q_n) \cup E$ again, which contradicts Corollary~\ref{cortang}.

Thus $C$ and $D$ are meromorphic at $0$. We have 
$$F^{-1}dF=\pmatrix{g&-g^2 \cr 1&-g \cr}\omega$$ 
so that $dC=(Cg+D)\omega$ and $dD=-g(Cg+D)\omega$. Consequently 
$g=-{dD \over dC}$ is also meromorphic at the puncture and 
this proves $E$ has finite total curvature. \phantom{almost}
\hfill\qed

\specialnumber{7} \proclaim{Lemma}
\label{periode}
Let $A${\rm ,} $B${\rm ,} $C${\rm ,} $D$ be the holomorphic {\rm (}\/multivalued\/{\rm )} data on $D^*$ 
parametrizing the end $E$ of Theorem~{\rm \ref{the3}.} Then 
$A(z)=z^{\beta}f(z)$ for some real $\beta$ and $f$ holomorphic on 
$D^*${\rm .} Also{\rm ,} $B,C$ and $D$ have similar representations{\rm .}
\endproclaim

{\it Proof}.
Let $\widetilde{D^*}=\{y \in \C;{\rm Re}\ y \leq 0\}$ and $e^y=z \in D^*_{\phantom{|}}$ be
the covering map. We have $F(y)=\pmatrix{A(y)&B(y) \cr C(y)&D(y) \cr}$ 
in ${\rm SL}(2,\C)$, and $F(y+2\pi i)=F(y)\,H$, where $H \in {\rm SU}(2)$ by 
Section~\ref{repres}.

Let $P \in {\rm SU}(2)$ diagonalize $H$, $PHP^{-1}=\triangle=
\pmatrix{e^{i \beta 2\pi}&0 \cr 0& e^{-i \beta 2\pi} \cr}$. Then 
$F_1=FP^{-1}$ defines the same end $E$ and 
$$F_1(y+2\pi i)=F(y+2\pi i)P^{-1}=F(y)HP^{-1}=F(y)P^{-1}\triangle=
F_1(y)\triangle.$$

Thus $A_1(y+2\pi i)=A_1(y)e^{i \beta 2\pi}$ and similarly for $B_1$, 
$C_1$ and $D_1$.

Now define $f(y)=e^{-y\beta}A_1(y)$, so that 
$$f(y+2\pi i)=e^{-(y+2\pi i)\beta}A_1(y+2\pi i)=f(y),$$ and $f$ defines 
a holomorphic map $f(z)$ on $D^*$, by $f(z)=f(y)$, $e^y=z$.

Then $e^{y \beta}f(z)=A_1(y)$, so that the (multi-valued) $A_1(z)$ on $D^*$ 
satisfies $A_1(z)=z^{\beta}f(z)$.
\hfill\qed

\specialnumber{4} \proclaim{Theorem}
\label{the4}
Let $E$ be a properly embedded Bryant annular end{\rm .} If $\b_{\infty}E$ is 
not $S_{\infty}$ then $E$ has finite total curvature{\rm .}
\endproclaim

\specialnumber{3} \proclaim{{C}orollary}
\label{standardends}
Let $E$ be a properly embedded Bryant annular end{\rm .} If $\Enod${\rm ,} then 
$E$ is regular and the total curvature of $E$ is finite{\rm .} $E$ is 
asymptotic to a catenoid cousin or horosphere end{\rm .}
\endproclaim

{\it Proof of the corollary}.
Theorem~\ref{the4} and Corollary~\ref{cororeg} yield the facts that $E$ is regular, 
$\b_{\infty}E$ is one point and the total curvature of $E$ is finite. 
Then the theorem of E.~Toubiana and R.~Sa~Earp yields the asymptotic 
behavior \cite{Tou-Ric}. \phantom{swim}
\hfill\qed

\demo{Proof of the theorem}
We know from Theorem~\ref{the3}, that $E$ will have finite total curvature 
if we can find a catenoid cousin $C_\alpha$ with $E$ on the mean convex side 
of $C_\alpha$; we will find such a $C_\alpha$ to prove Theorem~\ref{the4}. 
By Theorems~\ref{the1} and~\ref{the2} we know that $\b_{\infty}E$ is one 
point, which we take to be infinity in the upper half-space model of $\h$.

Let $B$ be a ball in $\h$, whose interior contains $\b E$ and where $E$ is 
transverse to $\b B$. Let $E_1$ denote the noncompact component of $E-B$, 
and let $W$ denote the mean convex domain (along $E_1$) bounded by $E_1$ 
and a compact domain on $\b B$.

If $x_3 \geq c >0$ on $E$ then $C_\alpha$ can be constructed using a 
catenoid cousin end below height $c$ which is a graph over an exterior 
domain $x_1^2+x_2^2>r_0^2$, asymptotic to the plane $x_3=0$ at infinity. 
So we can assume there is a sequence $q_n \in E_1$ with 
$x_3(q_n) \rightarrow 0$. Since $\b_{\infty}E=\infty$, we have 
$r(q_n)=\sqrt{x_1(q_n)^2+x_2(q_n)^2} \rightarrow \infty$.

For $q \in E_1$, let $\gamma$ be the minimizing geodesic of $\h$ joining 
$q$ to a point of $\b B$. We will be working with $q$ lower than $B$. 
Assume  $B=\{x_1^2+x_2^2 +\break(x_3-4)^2=1\}$ for convenience, and 
$x_3(q) \leq 1$, $r(q)>6$. Parametrize $\gamma$ by arc length so that 
$\gamma(0)$ is the highest point of $\gamma$ (which is not on $B$ by our 
choice of constants), and $\gamma(t_0)=q$ with $t_0<0$.

Let $P(t)$ be the family of (hyperbolic) planes orthogonal to $\gamma$ at 
$\gamma(t)$. For $t$ very negative, $P(t)$ is disjoint from $E_1$ since 
$\b_{\infty}E_1=\infty$, and $E_1$ is proper so that there is a first 
$t_1 \leq t_0$ (as $t$ increases) such that $P(t_1)$ touches $E_1$ at a 
point $q_1$.

We do Alexandrov reflection of $E_1$ with the planes $P(t)$ as $t$ 
increases from $t_1$ to $0$. Let $S(t)$ be symmetry of $\h$ through 
$P(t)$, $E_1(t)^+$ the part of $E_1$ on the side of $P(t)$ not containing 
$B$, and $E_1(t)^*=S(t)\left( E_1(t)^+ \right)$.

For $t$ slightly larger than $t_1$, $E_1(t)^+$ is a graph over (part of) 
$P(t)$,\break ${\rm int}\left(E_1(t)^*\right) \subset W$, and the angle between $P(t)$ 
and $E_1(t)^+$ is never $\pi/2$ along $\b E_1(t)^+$. These properties 
continue to hold until the first $t$ ($t_2$ say) such that $E_1(t_2)^*$ 
touches $\b B$, for if one of these properties failed to hold at some 
earlier~$t$, $P(t)$ would be a plane of symmetry of $E$. Then $E$ is part 
of a properly embedded, mean curvature one, compact surface $M$, with 
$\b M=\emptyset$. This is impossible.

Clearly $t_2<0$ since $q$ is lower than $B$, and so the symmetry of $q$ 
through some plane $P(t)$, $t<0$ meets $B$. Thus there is some point 
$\widetilde{q} \in E_1(t_2)^+$ such that $S_{t_2}(\widetilde{q})\in B$.

Let $\delta_1={\rm dist}(\widetilde{q},\gamma)$, and $q_t=S_t(\widetilde{q})$. 
Since $\gamma$ is invariant by $S_t$, we have ${\rm dist}(q_t,\gamma)=\delta_1$ 
as well. For $t=t_2$, $q_t$ is on $\b B$, so that
${\rm dist}(q_t,\gamma) \leq {\rm diam}(B)=\delta$. The curve $q_t$ joining 
$\widetilde{q}$ to $\b B$, as $t$ varies from $t_1$ to $t_2$, is an 
equidistant curve~$\beta$ whose distance from $\gamma$ is less than 
$\delta$, and this equidistant curve is contained in~$W$. We emphasize that 
this discussion is valid for any $q \in E_1$ with $x_3(q)<1$, $r(q)>6$.

In particular, consider the sequence $q_n \in E_1$, satisfying 
$x_3(q_n) \rightarrow 0$, $r(q_n) \rightarrow \infty$. Then a subsequence 
of the geodesics $\gamma_n$ joining $q_n$ to $B$ converges to a vertical 
geodesic over $B$ and the equidistant curves $\beta_n $ from 
$\widetilde{q}_n$ to $B$ are in $W$ and a distance at most $\delta$ from 
$\gamma_n$. So the equidistant curves $\beta_n$ are in the tubular 
neighborhood of $\gamma_n$ of radius $\delta$. As $n \rightarrow \infty$, 
the tubular neighborhoods converge to a vertical cone of hyperbolic width 
$\delta$. Let $C(\delta)$ denote this cone; for simplicity we can assume 
the base of $C(\delta)$ is the origin.

Now we can prove that $E_1 \cap A$ is a graph where 
$A=\left\{ x_3<1,r\geq6 \right\}$.

Suppose this were not true. Let $N$ be the Euclidean unit normal to $E$, 
$N.\veh>0$ and suppose that $N_3 \leq 0$ at some point $q \in E_1 \cap A$. 
Then the horosphere tangent to $E$ at $q$, $H(q)$, is at most of 
(Euclidean) radius $1$ and $\b E \subset H(q)^-$.

Then by Corollary~\ref{cortang}, $H(q)$ separates $W$ into three connected 
components. One is compact and contains part of $\b B$. One is noncompact, 
and contains the points $\widetilde{q}_n$, $n$ large. And the third is 
compact and inside $H(q)^+$. But the equidistant curves $\beta_n$ are 
in $W$ and disjoint from $H(q)$ for $n$ large; this is impossible since 
the $\beta_n$ go to $\b B$ in $W$. This proves $E_1 \cap A$ is a graph.

In fact the above argument proves much more: for $q \in A$, $H(q)$ must 
intersect $C(\delta)$; otherwise the equidistant curves $\beta_n$ would be 
disjoint from $H(q)$ for $n$ large; cf.\ Figure~9.

For $q \in E_1 \cap A$, let $R$ be the Euclidean radius of $H(q)$ and let 
$d$ be the Euclidean distance of $q$ to $C(\delta)$. Then (since $C(\delta)$ 
is invariant by homothety from $\sigma$ and 
$C(\delta) \cap H(q) \neq \emptyset$) there is a $\lambda>0$ such that 
$$2R \geq d \geq 2 \lambda r(q),$$
and $\lambda$ depends only on $C(\delta)$. In particular 
$R \rightarrow \infty$ when $r(q) \rightarrow \infty$.

Now we shall prove that $E$ is below some horosphere $x_3={\rm constant}$.

We know that $E_1 \cap A$ is the graph of a function $u$ and in 
Theorem~\ref{the2}, we derived the formula:
$${\vert \nabla u \vert^2 \over u(1+\vert \nabla u \vert^2)}=
{2R-u \over R^2} \leq {2 \over R}.$$

\begin{center} \BoxedEPSF{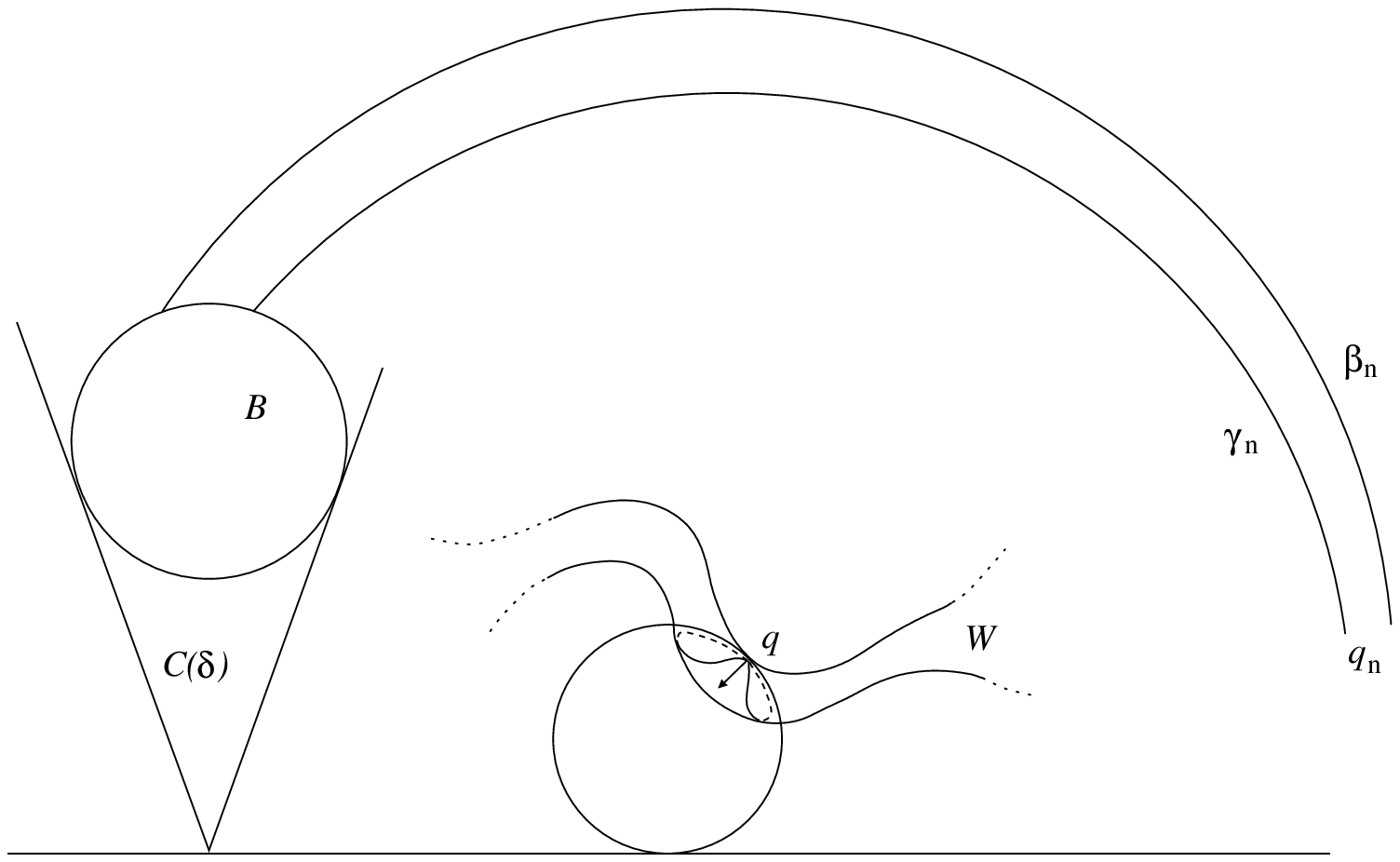 scaled 600} \end{center}
\centerline{Figure~9} 
\medbreak

Since $u \leq 1$ this implies 
$$\vert \nabla u \vert^2 \leq {2u \over R-2} \leq 
{2u \over \lambda r(q)-2}.$$
In particular, at the point $q_n \in E_1$, where 
$x_3(q_n) \rightarrow 0$, $r(q_n) \rightarrow \infty$, we obtain 
$$\vert \nabla u(q_n) \vert \leq \varepsilon^2_n,$$
for a sequence $\varepsilon_n \rightarrow 0$.

Now recall our discussion of Alexandrov reflection by planes orthogonal 
to the geodesics $\gamma_n$ joining $q_n$ to $\b B$. We found a point 
$\widetilde{q}_n$ in $E_1$, associated to the first accident of Alexandrov 
reflection, and we showed the equidistant curve $\beta_n$ from 
$\widetilde{q_n}$ to $\b B$ was in $W$. We have 
$\left| r(\widetilde{q}_n)-r(q_n) \right|<1$ by construction, so at 
$\widetilde{q}_n$ we also have an estimate 
$$\vert \nabla u(\widetilde{q}_n) \vert \leq \varepsilon^2_n,$$
for $\varepsilon_n \rightarrow 0$, 
$\varepsilon_n \sim {1 \over r(\widetilde{q}_n)^{1/4}}.$

Then the maximum oscillation of $u$ on the horizontal (Euclidean) disk $D$ 
of radius ${x_3(\widetilde{q}_n) \over \varepsilon_n}$, centered at 
$\widetilde{q}_n$, is $2\varepsilon_n x_3(\widetilde{q}_n)$.

To check this, notice that the most $\vert \nabla u \vert$ can be is 
$r(\overline{q})^{-1/2}$, where $\overline{q}$ is a point of $D$ closest 
to the origin. Thus,
$$r(\overline{q})=r(\widetilde{q}_n)-{x_3(\widetilde{q}_n) \over 
\varepsilon_n} \geq {r(\widetilde{q}_n) \over 2}.$$

Then $r(\overline{q})^{-1/2} \leq \sqrt{2}r(\widetilde{q}_n)^{-1/2}$ and 
the oscillation on $D$ is at most 
\[ \begin{array}{rl}
\displaystyle{\vert \nabla u(\overline{q}) \vert 
{x_3(\widetilde{q}_n) \over \varepsilon_n}} &\leq 
\displaystyle{\sqrt{2}r(\widetilde{q}_n)^{-1/2}x_3(\widetilde{q}_n) 
r(\widetilde{q}_n)^{1/4}} \\
&\leq 2x_3(\widetilde{q}_n)\varepsilon_n.
\end{array} \]

Define $D_n=D+(0,0,x_3(\widetilde{q}_n))$; $D_n$ is a horizontal disk 
above the graph of $u$ over $D$ so that $D_n \subset W$ and the hyperbolic 
radius of $D_n$ tends to infinity (it is $1/2\varepsilon_n$). Also 
the hyperbolic distance between $D_n$ and the graph of $u$ over $D$ 
is bounded by $\ln (2)$.

Let $t_n<0$ denote the first time that $S(t_n)(\widetilde{q}_n)$ touches 
$\b B$ (the first accident when we do Alexandrov reflection with the 
planes orthogonal to $\gamma_n$). We have $F_n=S(t_n)(D_n) \subset W$ and 
the distance of $F_n$ to $\b B$ is at most $\ln(2)$. As 
$n \rightarrow \infty$, $F_n$ converges to a horizontal horosphere $F$ 
which must be in $W$. Thus $E$ is below $F$.

Next we observe that $E_2=E \cap (\Omega\times\R^+)$ is a vertical graph, 
where $\Omega=\left\{ x_1^2+x_2^2>a^2 \right\}$, for some $a>0$. To see 
this, remark that $x_3(q) \leq c_0$ for some constant $c_0$ and so if 
$\veh(q)$ does not point up then $q$ is in the upper hemisphere of its 
tangent horosphere so $x_3(q) \geq R$ = the Euclidean radius of $H(q)$. 
Hence $R \leq c_0$ and $H(q)$ will be disjoint from the cone $C(\delta)$ 
for $r(q)$ larger than some fixed $a$. As before, this is impossible since 
the equidistant curves $\beta_n$, for $n$ large, will not intersect $H(q)$.

Now on the domain $\Omega \times \R^+$ where the subend $E_2$ is a graph, 
we consider the family of catenoid cousin ends $C(t)$ with each $C(t)$ a 
graph over $\Omega \times \R^+$, tangent to the vertical cylinder 
$\b \Omega \times \R^+$ and $\b C(t)$ is at height $t$ on 
$\b \Omega \times R^+$. These surfaces are described in \cite{Go}.

For $t>c_0$, $\b C(t)$ is above $E_2$. If $C(t)$ intersects $E_2$, then by 
Theorem~\ref{the5}, $\Gamma=C(t) \cap E_2$ is compact. Note that $\Gamma$ is not 
homologous to zero on $E_2$ (nor is any subcycle of $\Gamma$) since this 
would yield a compact domain $N$ on $E_2$ whose boundary is in $C(t)$. Now 
vary $t$ to obtain a last point of contact of $C(t)$ with $N$; then 
$C(t)=E_2$ by the maximum principle. It follows that $\Gamma$ is a Jordan 
curve on $E_2$ that generates $\Pi_1(E_2)$. On $C(t)$, $\Gamma$ bounds a 
catenoid cousin end that is below $E_2$ and Theorem~\ref{the4} is clear 
by Theorem~\ref{the3}.

Now, we can assume $C(t) \cap E_2=\emptyset$ for $t>c_0$, and then decrease $t$ 
to $0$. There is some largest $t$ where $C(t)$ is disjoint from $E_2$ and 
$C(s) \cap E_2 \neq \emptyset$, for $s<t$. Since $C(t)$ is vertical along 
$\b \Omega \times \R^+$ and $E_2$ is a graph (not vertical) there, 
$\b C(t)$ is always above $E_2$. Thus we are in the previous situation 
where $C(s) \cap E_2 \neq \emptyset$ and $\b C(s)$ is above $E_2$ and 
Theorem~\ref{the4} is proved.
\enddemo

\specialnumber{5} \proclaim{Theorem}
\label{the5}
Let $\Omega$ be a noncompact domain in the plane $\left(x_1,x_2\right)$ 
with at least one component of $\b \Omega$ noncompact{\rm .} Let $u_1,u_2$ be 
defined on $\Omega$ with their graphs solutions of the mean curvature 
equation $H=1$ in $\h${\rm .} Suppose the following conditions are satisfied\/{\rm :}

\begin{description}
\item[{\rm a)}] $u_2 \leq u_1 \leq 1$ on $\Omega${\rm ,} $u_1=u_2$ on $\b 
\Omega${\rm ,}

\item[{\rm b)}] $\displaystyle{{C_1 \over r^{\alpha}} \leq u_2 \leq 
{C_2 \over r^{\a}}}${\rm ,} for some positive constants $C_1,C_2,\a$ ($u_2$ 
is the graph of a catenoid cousin\/{\rm ),}

\item[{\rm c)}] $\displaystyle{{\vert \nabla u_1 \vert^2 \over u_1}
\leq {C \over r^2}}$, for some $C>0${\rm ,} $r^2=x_1^2+x_2^2${\rm .}
\end{description}

It then follows that $u_1=u_2$ on $\Omega${\rm .}
\endproclaim 

\demo{{R}emark {\rm 3}}
In order to apply this theorem to prove Theorem~\ref{the4}, we need to 
verify that the graph $u$ $(=u_1)$ of $E_2$   in Theorem~\ref{the4} 
satisfies the conditions a, b, and c. The conditions a and b are satisfied 
by construction; the condition c needs some discussion.
\enddemo

In the proof of Theorem~\ref{the2} we derived the gradient bound for $u$:
$${\vert \nabla u \vert^2 \over u\left(1+\vert \nabla u \vert^2 \right)} 
\leq {2 \over R},$$
where $R$ is the Euclidean radius of the horosphere $H(q)$. Since 
$u \leq 1$,  
$${\vert \nabla u \vert^2 \over u} \leq {2 \over R-2}.$$
So we need to know $R$ is of order $r^2$ for the graph $u$, to satisfy 
condition c.

We see this by considering $H(q)$, $q$ on the graph of $u$. Let $E$ denote 
the graph of $u$ (this is the $E_2$ in the proof of Theorem~\ref{the4}), 
and let $C$ be the vertical compact cylinder joining $\b E$ to the plane 
$x_3=0$. Observe that for $q \in E$, $H(q)$ must intersect $C$. For if 
$H(q)$ passes over $C$, then $\b E \subset H(q)^-$ so the figure eight in 
$H(q) \cap E$, would contain a Jordan curve $C_1$ that is homologous to 
$\b E$ on $E$ (Proposition~\ref{protang}). However, $u$ takes its maximum 
value on $\b E$ ($u$ has no interior maximum since the graph of $u$ would 
touch a horizontal horosphere at a local maximum and have the same mean 
curvature vector). Thus $C_1$ would be lower than $\b E$ and link the 
cylinder $C$.  Hence $H_q$ must intersect $C$. We want to estimate $1/R$ 
from above, so that for $q \in E$, we can assume $H(q)$ intersects the 
vertical segment over the origin at a point $p$ at height $x_3(p)$ less 
than some fixed $b>0$. Now for $x_3(q)<b$, the horosphere $H(q)$ passing 
through $p$, intersects the plane $x_3=0$ at the point $G(q)$; 
see Figure~10.

\begin{center} \BoxedEPSF{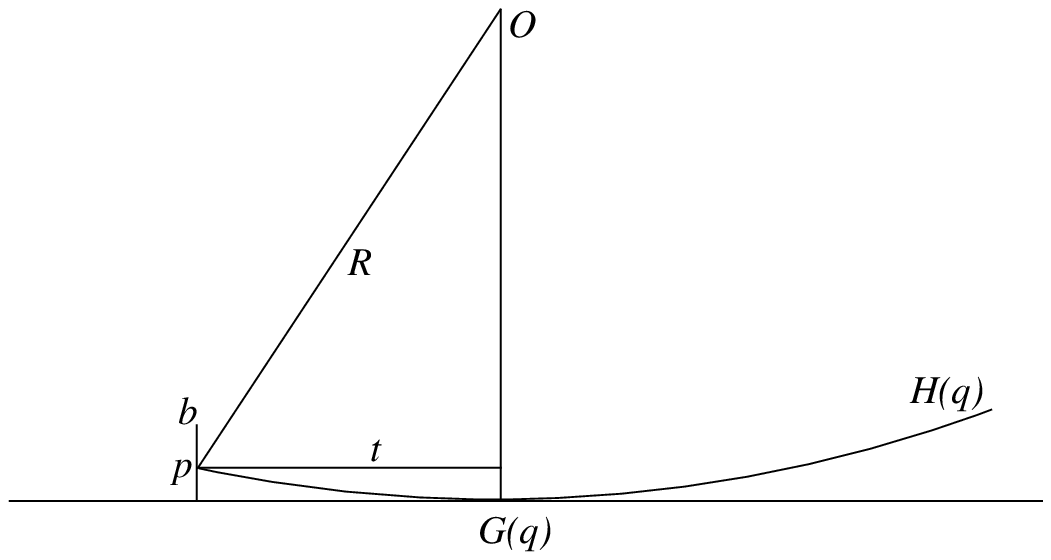 scaled 750} \end{center}
\centerline{Figure~10} 
\medbreak

Then $t^2+(R-x_3(p))^2=R^2$, so that 
$R={t^2 \over 2x_3(p)}+{x_3(p) \over 2} \geq {t^2 \over 2b}$. 
Apply the same (Pythagorean) calculation with $p$ replaced by $q$ to obtain 
$R \geq {\tau^2 \over 2b}$ where $\tau$ is the horizontal distance from $q$ 
to $G(q)$. Since $\tau+t \geq r$, $\tau$ or $t$ is at least $r \over 2$ 
so that $R \geq {r^2 \over 8b}$ as desired.

Before proving Theorem~\ref{the5}, we need some lemmas.

\specialnumber{8} \proclaim{Lemma}
\label{divergence}
Let $u$ be a solution of the equation $H=1$ on $\Omega${\rm .} Then $v=\ln u$ 
satisfies\/{\rm :}
$${\rm div} \left( {\nabla v \over W} \right)={-\vert \nabla u \vert^4
\over u^2W \left( 1+W \right)^2}\hbox{, where }
W^2=1+\vert \nabla u \vert^2.$$
\endproclaim

\demo{Proof}
We have $\displaystyle{{\rm div} \left( \nabla u \over W \right)={2 \over u}
\left(1-{1 \over W} \right)}$, ($H=1$), hence 
\[ \begin{array}{ll}
\displaystyle{{\rm div} \left( \frac{\n v}{W} \right)}&
\displaystyle{={\rm div} \left( \frac{\n u}{uW}\right)=
\frac{-\nor{u}{2}}{u^2W}+\frac{2}{u^2}\left( 1-\frac{1}{W} \right)} \cr
& \cr
& \displaystyle{=\frac{-\nor{u}{2}}{u^2 W}+\frac{2}{u^2}\frac{\nor{u}{2}}
{W \left( 1+W \right)}=\frac{\nor{u}{2}}{u^2W}
\left[ \frac{2}{1+W}-1 \right]}\cr
& \cr
& \displaystyle{=\frac{\nor{u}{2}}{u^2W} \left( \frac{1-W}{1+W} \right)=
\frac{\nor{u}{2}}{u^2W} \left[ \frac{-\nor{u}{2}}{\left( 1+W \right)^2} 
\right]=\frac{-\nor{u}{4}}{u^2W \left( 1+W \right)^2}.}
\cr
\noalign{\vskip-24pt}
\end{array} \]
\enddemo
 
\specialnumber{9} \proclaim{Lemma}
\label{croissance}
Let $\Omega(r)=\left\{ x \in \Omega;\vert x \vert \leq r \right\}$ and 
$C(r)=\Omega(r) \cap \left\{ \vert x \vert=r \right\}${\rm .} Define 
$v=\ln u_1-\ln u_2$ {\rm  (}$u_1,u_2$ as in Theorem~{\rm \ref{the5}),} and 
$M(r)=\sup \left\{ \vert v(x) \vert;\vert x \vert=r \right\}${\rm .} Then if 
$v \neq 0$ there is a $\beta<\a$ such that $M(r) \geq (\a-\beta)\ln r${\rm .}
\endproclaim

\demo{Proof}
Consider a family of catenoid cousin graphs $u_{\tau}(x)$, with $u_\tau$ strictly above $u_2$ on $\partial\Omega$ and
$u_\tau$ comes down to $u_2$ as $\tau\to\alpha$, with $u_\tau=u_2$ for $\tau=\alpha$.   Parametrize so that the growth
of
$u_{\tau}$ is 
$1/r^{\tau}$, $\tau<\a$. As $\tau \rightarrow \a$, one cannot have 
$u_{\tau}$ above $u_1$ for all $\tau$  (otherwise    $u_1=u_2$).  Hence  the graph 
of some $u_{\beta}$, $\beta<\a$, intersects the graph of $u_1$. As usual, 
we know the intersection cannot be homologous to zero on the graph 
(vary $\tau$ to get a last point of contact), and the intersection is not 
one compact cycle (otherwise there is a catenoid cousin below $u_1$ and 
Lemma~\ref{croissance} is proved) so that  the intersection is not compact and 
$u_1$ is above $u_{\beta}$ on a noncompact domain. Thus 
$M(r) \geq \left( \a-\beta \right)\ln r$.
\enddemo

\demo{Proof of   Theorem {\rm 5}}
We study $v=\ln u_1-\ln u_2=v_1-v_2$. Clearly $v \geq 0$, $v=0$ on 
$\b \Omega$ and $v \leq \gamma \ln r$ for some positive $\gamma$. We will 
show that if $v$ is not identically zero, then for some integer $k>1$, 
$M(r)$ grows faster than $\left( \ln r \right)^k$. This latter growth is 
impossible and so $v \equiv 0$.

By Stokes' theorem, 
$$\int_{\Omega(R)} {\rm div} \left( v\frac{\n v_1}{W_1} \right)-
{\rm div} \left( v\frac{\n v_2}{W_2} \right)=\int_{\b \Omega (R)} 
v \langle \frac{\n v_1}{W_1}-\frac{\n v_2}{W_2}, N \rangle,$$
where $N$ is the outer conormal along $\b \Omega (R)$.
 Apply this equation to $v=v_1-v_2$, 
\begin{eqnarray}
&&\int_{\O(R)} \left( \n v_1-\n v_2 \right) \left( \frac{\n v_1}{W_1}-
\frac{\n v_2}{W_2} \right)+\int_{\O(R)} v\ {\rm div} \left( \frac{\n v_1}{W_1} 
\right)-v\ {\rm div} \left( \frac{\n v_2}{W_2} \right)
\label{eq:equation 1}\\
&&\hskip1.25in =\int_{\b\O(R)} v\langle \frac{\n v_1}{W_1}-\frac{\n v_2}{W_2},N\rangle.\nonumber 
\end{eqnarray}
By   Lemma~\ref{divergence}, and the estimates 
$\displaystyle{\frac{\nor{u_i}{4}}{u_i^2} \leq \frac{C_i}{r^4}}$, $i=1,2$, 
and $v \leq \gamma \ln r$, 
\begin{eqnarray*}
\left| \int_{R_0}^R \int_{C(r)} v\ {\rm div} \left( \frac{\n v_i}{W_i} \right) 
\right| &\leq& \left| \int_{R_0}^R \int_{C(r)} v \frac{\nor{u_i}{4}}
{u_i^2W_i \left( 1+W_i \right)^2} \right| \\ \noalign{\vskip4pt}
&\leq &
C_i \left| \int_{R_0}^R \int_{C(r)} \frac{\ln r}{r^4}\frac{1}
{W_i \left( 1+W_i \right)^2} \right|.
\end{eqnarray*}
Since the last integral converges we have 
$\displaystyle{\left| \int_{\O(R)} v\ {\rm div} \left( \frac{\n v_i}{W_i} \right) 
\right| \leq a_i}$, for some constants $a_1$, $a_2$. Then equation 
(\ref{eq:equation 1}) yields 
\begin{equation}
\label{eq:equation 2}
a_3+\int_{\O(R)} \left( \n v_1-\n v_2 \right) \left( \frac{\n v_1}{W_1}-
\frac{\n v_2}{W_2} \right) \leq \int_{\b\O(R)} v \langle 
\frac{\n v_1}{W_1}-\frac{\n v_2}{W_2},N \rangle.
\end{equation}

For $R_1>0$, define $\displaystyle{\mu(R_1)=\int_{\O(R_1)} 
\left( \n v_1-\n v_2 \right) 
\left( \frac{\n v_1}{W_1}-\frac{\n v_2}{W_2} \right)}.$
We have 
\begin{eqnarray*}
\left( \n v_1-\n v_2 \right) 
\left( \frac{\n v_1}{W_1}-\frac{\n v_2}{W_2} \right)&=&
W_1 \left( \frac{\n v_1}{W_1}-\frac{\n v_2}{W_2} \right)^2\\
\noalign{\vskip4pt}
&&+
\ \left( W_1-W_2 \right) \frac{\n v_2}{W_2}
\left( \frac{\n v_1}{W_1}-\frac{\n v_2}{W_2} \right).
\end{eqnarray*}
Also 
\[ \begin{array}{ll}
\displaystyle{\left| \left( W_1-W_2 \right) \frac{\n v_2}{W_2}
\left( \frac{\n v_1}{W_1}-\frac{\n v_2}{W_2} \right) \right|} &
\displaystyle{=\left| \frac{\n u_2}{u_2 W_2} \right| 
\left| \frac{\nor{u_1}{2}-\nor{u_2}{2}}{W_1+W_2} \right|
\left| \frac{\n v_1}{W_1}-\frac{\n v_2}{W_2} \right|} \cr\noalign{\vskip5pt}
& \displaystyle{\leq \frac{c_1}{r^3} \left| 
\frac{\n v_1}{W_1}-\frac{\n v_2}{W_2} \right|.}
\end{array} \]
Then (\ref{eq:equation 2}) implies 
\begin{eqnarray}
&&a_3+\mu(R_1) + \displaystyle{ \int_{R_1}^R \int_{C(r)} W_1 
\left| \frac{\n v_1}{W_1}-\frac{\n v_2}{W_2} \right|^2-
\int^R_{R_1} \int_{C(r)} {c_1 \over r^3}
\left| \frac{\n v_1}{W_1}-\frac{\n v_2}{W_2} \right|} \\ \noalign{\vskip5pt}
&&\hskip1.25in \displaystyle{\leq\  \int_{C(R)} v \langle \frac{\n v_1}{W_1}-
\frac{\n v_2}{W_2},N \rangle.}\nonumber
\end{eqnarray}

Define $\displaystyle{\eta(r)=\int_{C(r)} \left| \frac{\n v_1}{W_1} 
-\frac{\n v_2}{W_2} \right|}$. Now,
$${\eta^2(r)\over 2\pi r} \leq \int_{C(r)} \left| \frac{\n v_1}{W_1}-
\frac{\n v_2}{W_2} \right|^2 \leq \int_{C(r)} W_1 \left| 
\frac{\n v_1}{W_1}-\frac{\n v_2}{W_2} \right|^2.$$ 
Inequality (3) then implies 
\begin{equation}
\displaystyle{a_3+\mu(R_1)+\int_{R_1}^R {\eta^2(r)\over 2\pi r}-
\int^R_{R_1} {c_1\eta(r) \over r^3} \leq M(R)\eta(R).}
\end{equation}
Now $\b\Omega$ is not compact, $v=0$ on $\b C(r)$ and $M(r)$ is the maximum 
of $v$ on $C(r)$ so  that
$$M(r) \leq \int_{C(r)} \left| \n v \right|.$$
Next 
$$\left| 
\frac{\n v_1}{W_1}-\frac{\n v_2}{W_2} \right| \geq {1 \over W_1} 
\left| \n v_1-\n v_2 \right|-\vert \n v_2 \vert \left|
\frac{1}{W_2}-\frac{1}{W_1} \right|$$
and $1/W_1 \geq c_2>0$, $\vert \n v_2 \vert \leq {\alpha \over r}$, 
$\left| \frac{1}{W_2}-\frac{1}{W_1} \right| \leq 2$, so  that
$$\eta(r) \geq c_2M(r)-4\pi\alpha.$$
By Lemma~\ref{croissance} we conclude $\eta(r) \rightarrow \infty$, 
as $r \rightarrow \infty$, unless $v \equiv 0$. Then there is 
a constant $c_3>0$ and $R_0 \geq 0$ such that for $r \geq R_0$,
$${\eta^2(r) \over 2\pi r}-{c_1\eta(r) \over r^3} \geq {c_3\eta^2(r)\over r}.$$

Thus (4) may be replaced by (5) for $R_1 \geq R_0$:
\begin{equation}
\displaystyle{a_3+\mu(R_1)+c_3 \int^R_{R_1} {\eta^2(r) \over r} \leq 
M(R) \eta(R).}\end{equation}

Now we will show that for $R_1$ greater than or equal to some (other) 
$R_0$, we have $\widetilde{\mu}(R_1)=a_3+\mu(R_1)>0$, for 
$$\widetilde{\mu}(R_1)=a_3+\int_{\Omega (R_1)} W_1 
\left| \frac{\n v_1}{W_1}-\frac{\n v_2}{W_2} \right|^2+
\int_{\Omega(R_1)} \left( W_1-W_2 \right) {\n v_2 \over W_2} 
\left( \frac{\n v_1}{W_1}-\frac{\n v_2}{W_2} \right).$$
The module of the second integral is at most 
$$\int_0^{R_1} {c_1 \over r^3} \eta(r),$$
and $W_1 \geq 1$ so that $\widetilde{\mu}(R_1) \geq a_3+\int_0^{R_1} 
\left( {\eta(r)^2\over 2\pi r}-{c_1\eta(r) \over r^3} \right)$, 
which diverges since $\eta(r) \rightarrow \infty$.

Now $\widetilde{\mu}(R_1) \geq \widetilde{\mu} (R_0)+
c_3 \int^{R_1}_{R_0} {\eta^2(r) \over r}$. By Lemma~\ref{croissance} and 
the comparison between $\eta(r)$ and $M(r)$ we conclude 
$\widetilde{\mu}(R_1)$ grows at least as fast as $\ln^3(R_1)$.

We write equation (5) as:
\begin{equation}
\displaystyle{\widetilde{\mu}(R_1)+c_3 \int^R_{R_1} {\eta^2(r) \over r}
\leq A \eta(R),}
\end{equation}
for $R \in [R_1,R_2]$, and $A=\sup \{M(R);R_1 \leq R \leq R_2 \}$.

Let $\xi$ be the function defined on the interval $J=\left[ R_1,
R_1 \exp\left( {2A^2 \over c_3\widetilde{\mu}(R_1)}\right) \right)$ by 
$${c_3 \over A} \ln \left({R \over R_1}\right)={2A \over 
\widetilde{\mu}(R_1)}-{1 \over \xi(R)}.$$
On $J$, $\xi$ satisfies the equation:
$${\widetilde{\mu}(R_1) \over 2}+c_3 \int^R_{R_1} {\xi^2(r) \over r}=
A\xi(R).$$

The connected component of $\{R \in J \cap [R_1,R_2];\xi(R)<\eta(R) \}$ 
that contains $R_1$, is open by construction and closed by equation (6). 
Thus it is the interval $J \cap [R_1,R_2]$. Since $\xi(r) \rightarrow \infty$ 
when $r$ converges ($r$ increasing) to the right end point of $J$, and 
$\eta$ is bounded on $[R_1,R_2]$, we conclude $R_2 \in J$. Thus 
$$R_2 \leq R_1 \exp\left( {2A^2 \over c_3\widetilde{\mu}(R_1)} \right).$$ 
Since $A \leq \gamma \ln(R_2)$, we have 
$$\left[ {c_3\widetilde{\mu}(R_1) \over 2} 
\ln\left( R_2 \over R_1 \right) \right]^{1 \over 2} \leq \gamma\ln(R_2),$$ 
for $R_0 \leq R_1 \leq R_2$. However this contradicts our estimate for the 
growth of $\widetilde{\mu}(R_1)$ (take $R_2=R_1^2$). This completes
the proof of Theorem~\ref{the5}.
\enddemo

\section{Nondensity at infinity  of finite topology surfaces} 

Let $M$ be a properly embedded Bryant surface with $\b M$ perhaps not 
compact. Assume a properly embedded surface $\S$ exists with $\b \S=\b M$ 
and $\S \cup M=\b W$ with $M$ mean convex along $W$. Let $P$ be a 
(hyperbolic) plane with $\h-P=P^+\cup P^-$, the connected components of 
the complement. Assume $\S \subset P^-$. Let $M^+=M \cap P^+$.

\specialnumber{6} \proclaim{Theorem}
\label{the6}
There is a constant $c>0$ {\rm (}\/independent of $M$\/{\rm )} such that if 
$\vert K(q) \vert<c$ for $q \in M^+${\rm ,} then 
${\rm int}\left( \b_{\infty}M^+ \right)=\emptyset${\rm ;} i.e.{\rm ,} $M$ cannot be 
asymptotic to an open set at infinity in $P^+${\rm .} In the half-space model,
for $q \in M^+$ and $x_3(q)$ sufficiently small{\rm ,} $M^+$ is a vertical 
graph near $q${\rm ,} no point of $M$ is below this local graph{\rm ,} and the angle 
between $\veh(q)$ and $\overrightarrow{e_3}$ is at most $\pi / 4${\rm .}
\endproclaim

\demo{Proof}
We work in the upper half-space model. At each $q \in M$, $M$ is locally a 
graph over $H(q)$, in geodesic coordinates orthogonal to $H(q)$. If 
$\vert K \vert$ is small on $M$ then the second fundamental form of $M$ is 
close to that of $H(q)$, since $H=1$ and $K$ small implies the principal 
curvatures of $M$ are close to~$1$. Hence there is a $c>0$ such that if 
$q \in M$ and $\vert K(q) \vert<c$, then $M$ is a graph over the disk 
$D(q)$ of radius 3 in $H(q)$, centered at $q$, and the maximum distance 
of the graph to this disk $D(q)$ is one-half. We will see that this $c$ works 
in Theorem~\ref{the6}.

We now suppose $\vert K(q) \vert<c$ for $q \in M^+$. Let $q \in M^+$ 
and suppose $\veh(q).\overrightarrow{e_3} \leq 0$ (i.e., $\veh(q)$ 
points down). Then $H(q)$ is a Euclidean sphere tangent to $S_{\infty}$ 
at one point. The upper hemisphere of $H(q)$ has (hyperbolic) diameter 
2 and $q$ is in this upper hemisphere so that $D(q)$ contains this hemisphere. 
Hence $M$ is a graph over the upper hemisphere. We call this graph 
${\rm Cap}(q)$. The graph is at most a distance one-half from the hemisphere 
so that $x_3$ has a maximum at an interior point $p \in {\rm Cap}(q)$. At $p$, 
$\veh(p)$ has the direction of $-\overrightarrow{e_3}$ (by comparison 
with the horizontal horosphere $\{x_3=x_3(p)\}$) and a simple 
calculation of the Euclidean Gaussian curvature at a point of $M$ with 
$\veh$ parallel to $-\overrightarrow{e_3}$ shows $M$ is strictly 
Euclidean convex at $p$. So the planes $x_3={\rm constant}$ meet ${\rm Cap}(q)$ in 
convex compact curves at heights a little below $x_3(p)$.

We can assume $P^+=\left\{ x_1^2+x_2^2+x_3^2 \leq 9,x_3>0 \right\}$ and 
the origin $\sigma$ is in $\b_{\infty}M$. We will prove that if 
$q \in M^+$ and $q$ is sufficiently close to $\sigma$ (in the Euclidean 
metric) then $M$ is a vertical graph in a neighborhood of $q$, over a 
domain $\Omega \subset \left\{ x_3=0 \right\}$, and in 
$\Omega \times \R^+$, there is no point of $M$ below this graph of $M$ 
near $q$. Thus $\Omega \cap \b_{\infty}M=\emptyset$ and 
${\rm int}\left( \b_{\infty}M^+ \right)=\emptyset$.

Define ${\rm Cyl}(r)=\left\{ x_1^2+x_2^2<r^2,x_3>0 \right\}$ and suppose 
$q \in M^+\cap {\rm Cyl}(1/4)$. If $\veh(q). \overrightarrow{e_3}>0$ for each 
such $q$ with $x_3(q)$ sufficiently small then $M$ is a graph over a 
domain $\Omega$ and if $M_1$ denotes this part of $M$ near $q$ where $M$ 
is a vertical graph, then for $p \in M_1$, the vertical segment from $p$ 
to $\left\{ x_3=0 \right\}$ cannot meet $M$ again since at the first point 
where this segment again meets $M$, the vector $\veh$ would necessarily 
point into $W$, hence it would have to point down, a contradiction. Thus 
it suffices to prove $\veh(q).\overrightarrow{e_3}>0$ for $x_3(q)$ 
sufficiently small.

Suppose the contrary, $\veh(q).\overrightarrow{e_3} \leq 0$, for $q$ 
arbitrarily low. For $q \in {\rm Cyl}(1/4)$ and $x_3(q) \leq 1/8$, we know $H(q)$ 
has at most (Euclidean) radius $1/8$ so $H(q) \subset {\rm Cyl}(1)$, and 
${\rm Cap}(q) \subset {\rm Cyl}(1)$. Let $p \in {\rm Cap}(q)$ be a point where $x_3(p)$ is 
a local maximum and the level curves of ${\rm Cap}(q)$ near $p$ are compact 
Jordan curves $C(t)$ in the planes $x_3={\rm constant}$.

Consider the evolution of these level curves $C(t)$ as $x_3$ decreases 
from $x_3(p)$. For values  near $x_3(p)$, there is no other part of $M$ 
inside the disk $D(t)$ of $\left\{  x_3=t \right\}$ bounded by $C(t)$. 
As long as $C(t)$ stays compact and nonsingular, there is no other part 
of $M$ in $D(t)$, since the part would bound a compact domain above 
$x_3=t$, and under $\left\{ C(\tau);t \leq \tau \leq x_3(p) \right\}$ 
and at the highest point of this compact part of $M$, $\veh$ is parallel 
to $\overrightarrow{e_3}$ so that $M$ would equal a horosphere $x_3={\rm constant}$; 
cf.\ Figure~11-a.

Also notice that $C(t)$ cannot acquire a singularity (i.e., a point 
where $\nabla x_3=0$) as long as $C(t)$ stays compact. For if a singularity 
occurs at a point $q_1 \in C(t)$ then $\veh(q_1)$ is vertical. It cannot 
point up, since then $D(s)$ would contain other parts of $M$ for $s>t$, 
$s$ near $t$, and this is impossible by the previous paragraph; 
see Figure~11-b.
\figin{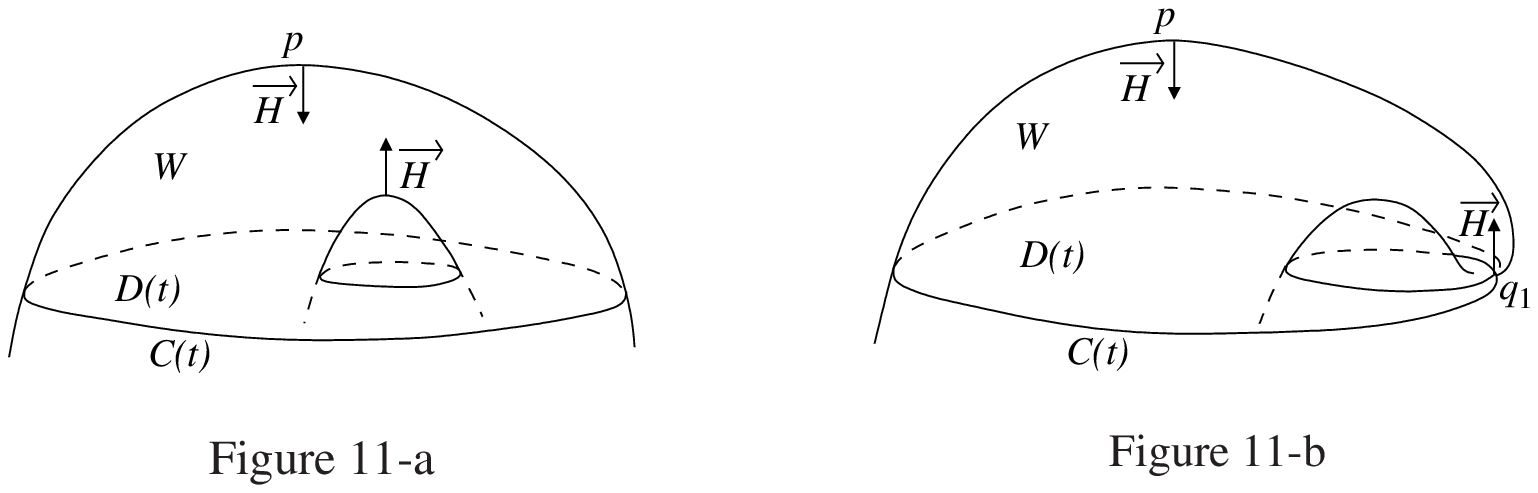}{780}

But $\veh(q_1)$ cannot point down either since $M$ would then be strictly 
locally (Euclidean) convex near $q_1$ and $x_3$ would have a local maximum 
at $q_1$, not a critical point of negative index.

Thus as long as $C(t)$ stays inside ${\rm Cyl}(2)$, it is a smooth Jordan 
curve.

As $t$ decreases to zero, $C(t)$ must leave ${\rm Cyl}(2)$ since ${\rm Cap}(q)$ must 
connect to the rest of $M$. Thus there are values of $t$ where $C(t)$ 
traverses $\b {\rm Cyl}(3/2)$.

Now if $x_3(q)$ is small and $q \in {\rm Cyl}(1/8)$, there will be points 
$\widetilde{q}$ of $C(t)$ in $\b {\rm Cyl}(3/2)$ where 
$\veh(\widetilde {q}).\overrightarrow{e_3} \leq 0$ and 
$x_3(\widetilde {q}) \leq 1/8$. Then ${\rm Cap}(\widetilde{q}) \subset 
{\rm Cyl}(2)-{\rm Cyl}(1)$ and ${\rm Cap}(\widetilde{q})$ has a local maximum 
of $x_3$ near $\widetilde{q}$. So the curves $C(t)$ are not connected 
before leaving ${\rm Cyl}(2)$, a contradiction.

It remains to obtain the gradient bound for the graph. For any horosphere 
of (Euclidean) radius $R$ in $\h$, the part of the horosphere where the 
mean curvature vector makes an angle greater than $\pi/4$ with 
$\overrightarrow{e_3}$  is of hyperbolic diameter at most $5$. So if 
$q \in M^+$ and $\vert K \vert$ is sufficiently small on $M^+$, then $M$ 
will be a graph over a geodesic disk in $H(q)$ that contains the northern 
hemisphere of $H(q)$, if the angle between $\veh(q)$ and 
$\overrightarrow{e_3}$ is greater than $\pi/4$. Now the same argument as 
before (with ${\rm Cap}(q)$ and the $C(t)$) leads to a contradiction. This 
proves Theorem~\ref{the6}.
\enddemo

\specialnumber{7} \proclaim{Theorem}
\label{the7}
Let $E$ be a properly embedded Bryant annular end{\rm .} If 
$\b_{\infty}E=S_{\infty}$ then there is a proper arc $\gamma$ on $E$ 
with $\b_{\infty} \left\{ \gamma(t);t \geq 0 \right\}=p_1${\rm ,} 
$\b_{\infty} \left\{ \gamma(t);t \leq 0 \right\}=p_2$ and $p_1 \neq p_2${\rm .}
\endproclaim

\demo{Proof}
On any subend of $E$, there must be points $q$, with 
$2x_3(q)<{\rm inf}\left( x_3\res_{\b E} \right)$ and 
$\veh(q).\overrightarrow{e_3} \leq 0$; otherwise the subend would be a 
graph near $x_3=0$, and so could not be dense at  infinity.

At such a point $q$, $\b E \subset H(q)^-$ since $H(q)$ is a sphere of 
Euclidean radius less than $2x_3(q)$. If the connected component $\S(q)$ 
of $q$ in $E \cap H(q)$ is not compact then there is an arc $\gamma_q$ on 
$H(q)$ in $\S(q)$ joining $q$ to $G(q)$; i.e., $\gamma_q$ is asymptotic 
to $G(q)$ at infinity. We will show next that such a $q$ can be found so 
that $\S(q)$ is not compact.

Suppose $\S(q)$ is compact. Then Proposition~\ref{protang} gives $E_1$ 
and a figure eight $C_1 \cup C_2$, $C_1 \cup C_2=\b E_1$, 
$E_1 \subset H(q)^+$, $E_1$ is compact, and $E_1$ separates $E$. Assume 
$C_1$ is the Jordan curve homologous to $\b E$ in $E-E_1$; 
cf.\ Figure~5.

Let $E'$ be the subannulus of $E$ bounded by $C_1$. We can assume there 
are points $q' \in E'$, $q' \neq q$, with $G(q')=G(q)$, for we can 
consider $\widetilde{q}$ near $q$ on $E$; if we could not find 
$\widetilde{q}'$ on $E'$ with $G(\widetilde{q}')=G(\widetilde{q})$ then 
$G\res_{E'}$ would miss an open set $\Omega$ (the open set being the image 
by $G$ of an open set about $q$ on $E$) in $S_{\infty}$, $\Omega$ a 
neighborhood of $G(q)$. However $E'$ is dense at infinity so there are 
points $y$ of $E'$ converging to $\Omega$ with 
$\veh(y).\overrightarrow{e_3} \leq 0$ (otherwise $E'$ would be a graph 
near $\Omega$), and then $G(y) \in \Omega$, for $x_3(y)$ small.

So we can assume there is $q' \in E'$, $q' \neq q$ and $G(q')=G(q)$. If 
$\S(q')$ is not compact then the arc $\gamma_{q'}$ joining $q'$ to $G(q')$ 
exists on $H(q')$. Thus, we suppose $\S(q')$ compact.

There are two possibilities.

\smallbreak
$\bullet$ Case 1. $\b E' \subset H(q')^+$. In this case, we have a compact 
component $E_1' \subset \left(E'\cap H(q')^+\right)$ with $\b E'_1=C_1'$. 
And there is a compact disk $D' \subset H(q')$ with $D' \cup E_1'=\b Q_1'$, 
$Q'_1$ a compact domain in $H(q')^+$ (Proposition~\ref{protang} and 
Figure~4). Also $\b E' \subset Q_1'$.

\smallbreak

Now $Q_1' \cap H(q)^+$ contains a connected compact component $Q$ with\break
$\b E' \subset Q$; cf.\ Figure~5. $Q$ is mean convex and 
$E_1 \subset Q$ so that $Q_1 \subset Q$. Also, $Q_1$ is  mean convex along $E_1$. 
Since $E_1$ and $E_2=\left( \b Q \right) \cap {\rm int}\left( H(q)^+ \right)$ 
are on $E$, there must be another component $F$ of $E$ in $Q-Q_1$ that 
separates $E_1$ and $E_2$ ($W$ is mean convex along $E$). Then $\veh$ 
points into the noncompact component of $H(q)^+-F$, along $F$, and 
this contradicts Lemma~\ref{lemins}.

\smallbreak
$\bullet$ Case 2. $\b E' \subset H(q')^-$. In this case, a Jordan curve 
of $H(q')$, $C'_1$ say, together with $C_2$ bounds a compact annulus 
$N \subset E$. Near $C_2$, $N$ is outside $H(q)$, so that $N \cap H(q)^-$ is 
a compact domain on $E$ with boundary on $H(q)$ and outside $H(q)$. This 
contradicts Lemma~\ref{lemout}. Thus we can construct a proper arc 
$\gamma_q$ from $q$ to $G(q)$ on $E$.
\smallbreak

Now do the same construction at a point $q_1 \in E$ with $G(q) \neq G(q_1)$. 
Join $q$ to $q_1$ by a path $\delta$ on $E$. Then the arc 
$\gamma=\gamma_q \cup \gamma_{q_1} \cup \delta$ works to prove the theorem.
\enddemo

\specialnumber{8} \proclaim{Theorem}
\label{the8}
Let $M$ be a properly embedded Bryant surface{\rm .} Suppose $\gamma$ is a 
proper arc on $M$ that separates $M$ into two components $M_1${\rm ,} $M_2${\rm .} There 
exist two properly embedded Bryant surfaces $\S_1${\rm ,} $\S_2$ satisfying\/{\rm :}
\begin{itemize}
\item[{\rm a)}] $\S_1$ and $\S_2$ are stable{\rm ,} $\b \S_1=\b \S_2=\gamma${\rm ,} 
$\S_1 \cap \S_2=\gamma${\rm ,} 
\item[{\rm b)}] $\S_1 \cup \S_2$ bounds a domain $R$ contained in the 
mean convex component $W$ of $\h-M${\rm ,} 
\item[{\rm c)}] $R$ is mean convex{\rm ,} 
\item[{\rm d)}] $\S_1 \cup M_1$ separates $\h$ and $\S_2 \cup M_2$ as well{\rm .}
\end{itemize}
\endproclaim
 
\demo{Proof} 
Fix a point $p \in \gamma$ and let $B=B_R$ denote the ball of $\h$ 
centered at $p$ of radius $R$. Let $M_R$ be the connected component of 
$M \cap B$ containing~$p$. The connected component of $\gamma \cap M_R$ 
containing $p$, separates $M_R$ into two components; denoted $M_1(R)$ and 
$M_2(R)$.

$M_R$ together with a compact domain on $\b B \cup M_0$, $M_0$ the part 
of $M-M_R$ in $B$, bound a mean convex domain $Q$; $Q \subset W$. The part 
of $\b Q$ on $\b B$ has mean curvature greater than $1$.

Let $D_1 \subset Q$ be  a least area embedded minimal surface with 
$\b D_1=\Gamma_1=\b M_1(R)$, and let $Q_1$ be the compact domain bounded 
by $D_1 \cup M_1(R)$. $D_1$ is a barrier for the Plateau problem so we can 
find a least area minimal surface $D_2 \subset Q-Q_1$ with 
$\b D_2=\Gamma_2=\b M_2(R)$. Let $Q_2$ be the compact domain bounded by 
$D_2 \cup M_2(R)$. We have $Q_1 \cup Q_2 \subset Q \subset W$ and 
${\rm int}(Q_1) \cap {\rm int}(Q_2)=\emptyset$.

Now consider domains $\widetilde{Q} \subset Q$ with 
$\b \widetilde{Q}=M_1(R) \cup \S$, $\S$ a surface with 
$\b \S=\b M_1(R)=\Gamma_1$. The functional on 
$\left(\widetilde{Q},\b \widetilde{Q}\right)$: 
$$(\widetilde{Q},\b \widetilde{Q}) \mapsto 
{\rm area}\left( \S\right)+2{\rm Vol}\left( \widetilde{Q} \right)$$
has a minimum and at such a $\widetilde{Q}$, the smooth points of $\S$ 
have mean curvature-one. This is proved in \cite{A-R} when the mean 
curvature of $\b Q$ is strictly greater than one; the only difference is 
that the minimum may now touch $M_1(R)$, in which case $\S=M_1(R)$ and 
$M_1(R)$ is stable in $Q$.

So let $\S_1$ be  a minimum, $\b \S_1=\Gamma_1$, $\S_1 \cup M_1(R)=
\b \widetilde{Q}$, and the mean curvature of $\S_1$ is one.

Observe that $\S_1 \subset Q_1$ (this is proved in \cite{A-R}) since, if 
$\widetilde{Q}$ went outside $Q_1$, one could remove the part of 
$\widetilde{Q}$ outside of $D_1$ and  reduce the functional.

Notice also that the mean curvature vector of $\S_1$ points outside of 
$\widetilde{Q}$. Otherwise $\widetilde{Q}$ would be mean convex so that one 
could find a least area minimal surface 
$\widetilde{D} \subset \widetilde{Q}$, $\b \widetilde{D}=\Gamma_1$. Then 
the  functional is smaller on the domain bounded by 
$\widetilde{D} \cup M_1(R)$; a contradiction.

Now working with $\Gamma_2=\b M_2(R)$ and $Q_2$, one finds a mean curvature-one, stable surface $\S_2 \subset Q_2$,
$\b \S_2=\Gamma_2$ and the mean  curvature vector of $\S_2$ points outside of the domain bounded by 
$\S_2 \cup M_2(R)$. Thus the domain of $W \cap B$ bounded by 
$\S_1 \cup \S_2$ (and a part of $\b B$) is mean convex; cf.\ Figure~12.

\figin{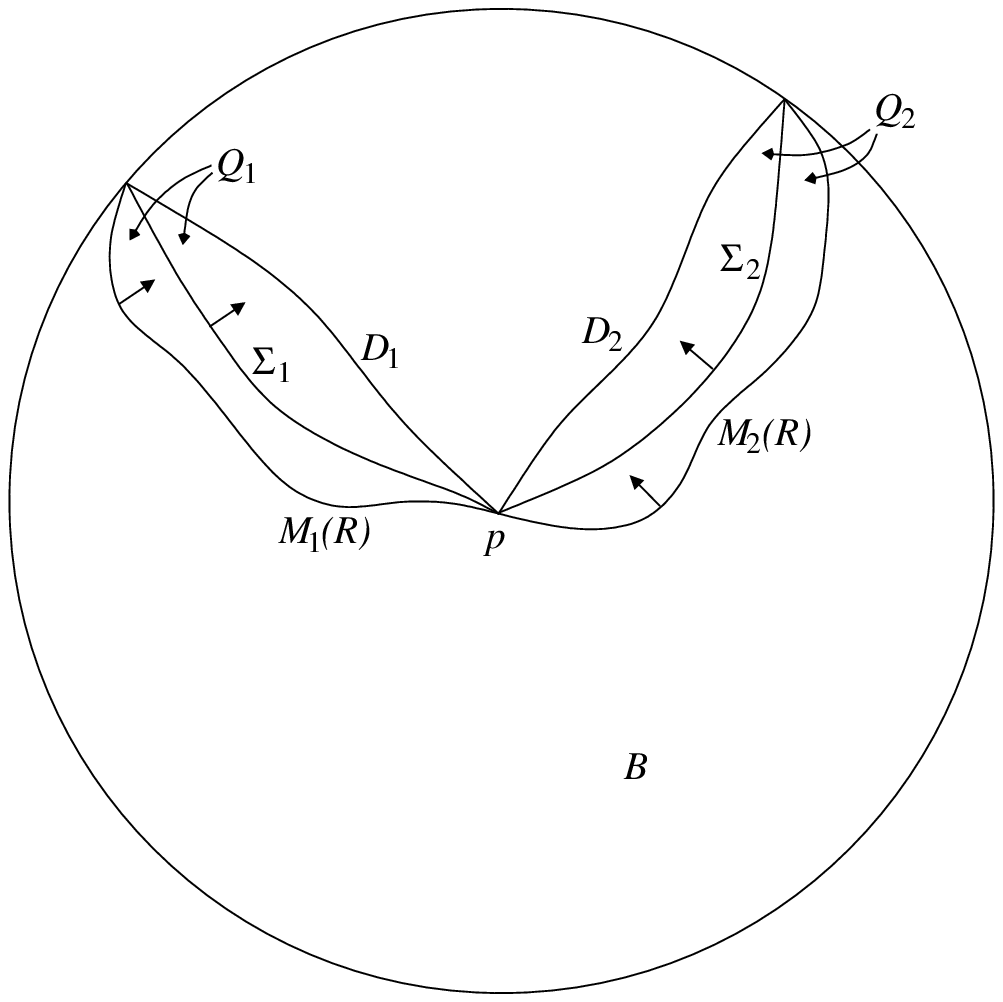}{700}
\centerline{Figure~12}
\medbreak

For $R>r>0$, one has uniform area and curvature bounds of $\S_1$ and 
$\S_2$ on balls of radius $r$ a fixed distance from $\b \S_1$ and 
$\b \S_2$. Then (as in \cite{A-R}), one can find a convergent subsequence 
of $\S_1$ and $\S_2$, as $R \rightarrow \infty$, which yield the $\S_1$ 
and $\S_2$ of Theorem~\ref{the8}.

In the case $\S_1 \mbox{ (or }\S_2) \subset M$ then this part of $M$ is 
stable in $W$ but this easily implies stability in $\h$ (look at an 
unstable domain $D$ corresponding to a first eigenvalue $\lambda_1<0$). 
This proves Theorem~\ref{the8}.
\enddemo

\specialnumber{9} \proclaim{Theorem}
\label{the9}
Let $M$ be a properly embedded Bryant surface of finite topology{\rm .} 
Then $\b_{\infty}M \neq S_{\infty}${\rm .}
\endproclaim

{\it Proof}.
Assume the contrary; $M$ is dense at infinity. Then by 
Corollary~\ref{cororeg} for some annular end $E$ of $M$, 
$\b_{\infty}E=S_{\infty}$, and so Theorem~\ref{the7} applies: there is a proper 
arc $\gamma$ on $E$ and $\b_{\infty}\gamma$ equals two distinct points 
$p_1$, $p_2$. Also, $E$ has genus zero so that  $\gamma$ separates $E$, hence $M$ as 
well. Theorem~\ref{the8} then yields stable surfaces $\S_1$, $\S_2$ 
satisfying the conditions a through d of Theorem~\ref{the8}.

Let $\Gamma \subset S_{\infty}=\left\{x_3=0 \right\} \cup 
\left\{ \infty \right\} $, be a circle separating $p_1$ and $p_2$. Note that $\S_1$ 
and $\S_2$ are stable so their curvature is small far from $\gamma$. 
In particular, when $c$ is the constant of Theorem~\ref{the6}, there is a 
$c_0>0$ such that $\vert K(q) \vert<c$ for $q \in \S=\S_1 \cup \S_2$, 
${\rm dist}(q,\gamma) \geq c_0$.

Let $\cal T$ be those points of $\h$ whose Euclidean distance to $\Gamma$ is 
at most $c_1>0$. Then for $c_1$ sufficiently small, $\S_1 \cap \cal T$ and 
$\S_2 \cap \cal T$ are vertical graphs over domains $\Omega_1$ and 
$\Omega_2 \subset \left\{ x_3=0 \right\}$. We know 
$\veh.\overrightarrow{e_3}>0$ on ${\cal G}=\left( \S_1 \cap \cal T \right) 
\cup \left( \S_2 \cap \cal T \right)$ and $\S_1 \cup \S_2$ bounds a mean 
convex domain $R$ by Theorem~\ref{the8}, so 
$\Omega_1 \cap \Omega_2=\emptyset$. Also we can assume the angle between 
$\veh(q)$ and $\overrightarrow{e_3}$ is less than $\pi/4$ on ${\cal G}$. 
Then for $q \in {\cal G}$, and $t=x_3(q)$ sufficiently  small, ${\cal G}$, 
near $q$, is a vertical graph over a horizontal disk $D(q)$, centered at 
$q$, of Euclidean radius~$t$.

Let $\tau>0$ and $\Gamma_\tau=\Gamma+\tau \overrightarrow{e_3}$. Choose 
$\tau$ small so that $ \Gamma_\tau \subset \cal T$ and $\Gamma_\tau$ is 
transverse to $\S$. The linking number of $\Gamma_\tau$ and $\gamma$ is 
one so that $\Gamma_\tau \cap \S_1$ consists of an odd number of points.  Now,  
$\S=\S_1 \cup \S_2$ bounds the mean convex domain $R$ so that there is an arc of 
$\Gamma_\tau$, which we denote $(q_1,q_2)$, joining a point $q_1 \in \S_1$ 
to $q_2 \in \S_2$ and the interior of the arc is in the interior of $R$.

For $q$ on the arc $(q_1,q_2)$, let $J(q)$ be the disk $D(q)$ together 
with the lower hemisphere of the horosphere that contains $\b D(q)$ and 
is vertical along $\b D(q)$. Note that $J(q)$ has a corner along $\b D(q)$.

For $q=q_2$, $J(q) \subset \Omega_2 \times \R^+$ by our gradient bound on 
the graph ${\cal G}$. Now move $q$ on the arc $(q_1, q_2)$ from $q_2$ to 
$q_1$. We know that $\Omega_1 \cap \Omega_2=\emptyset$ so that 
$J (q_2) \cap \S_1=\emptyset$. There will be a first $\widetilde{q}$ on 
the arc where $J(q)$ touches $\S_1$. We will next see that $J(q)$ 
touches $\S_1$ at infinity.

Suppose $J(q)$ first touches $\S_1$ at a smooth point $p$ on the horosphere 
in $J(q)$. The mean curvature vector of  the horosphere points up at $p$, 
and the mean curvature vector of $\S_1$ points up at $p$ too. So the vectors 
are equal and $\S_1$ is a horosphere. This is impossible because the proper 
arc $\gamma$ is on $\S_1$ and $\gamma$ has two points at infinity, $p_1$ 
and $p_2$; the horosphere has one point at infinity.

Next suppose the first point $p$ where $J(q)$ touches $\S_1$ is on 
$\b D(q)$. We know that the horizontal segment in $D(q)$, joining $p$ to 
$q$ (which we call $[p,q]$) meets $\S_1$ only at $p$. Also this segment 
does not meet $\S_2$ because our gradient bound implies 
$[p,q] \subset \Omega_1 \times \R^+$.

Thus the segment $[p,q]$ is contained in $R$. The (Euclidean) tangent 
plane to $\S_1$ at $p$ is a support plane of $J(q)$ and $\veh(p)$ points 
up at $p$. This contradicts the fact that $R$ is mean convex: $\veh(p)$ 
points into $R$, and $[p,q]$ $(\subset R)$    is on the other side of 
the tangent plane than $\veh(p)$.

Thus there is a point $q$ on the arc where $J(q)$ touches $\S_1$ for the 
first time at a point $q_{\infty} \in \Gamma$; see \ Figure~13.
\begin{center}
\BoxedEPSF{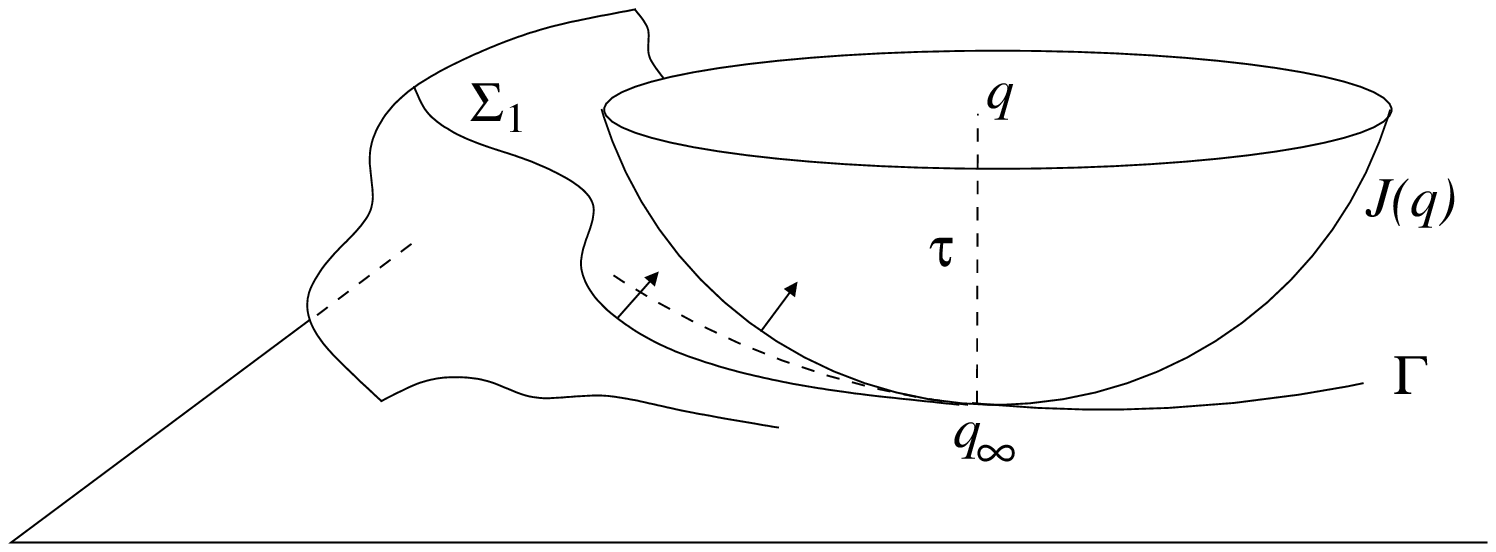 scaled 500} \end{center}
 \centerline{Figure~13} 
 \medbreak

Now consider $q'$ on the arc $(q,q_2)$ at Euclidean distance less than 
$\tau$ from $q$, such that the point of $\Gamma$ below $q'$ is not in 
$\b_{\infty}\Sigma_1$ but $q_\infty$ is below $D(q')$. By 
Lemma~\ref{catenoid}, there exists a one-parameter family of vertical 
graphs $C(t),\break 0<t \leq 1$, such that $C(1)$ is the original horosphere 
of $J(q')$, and $C(t)$ $(t<1)$ is a catenoid cousin end; each $C(t)$ is 
vertical along $\b C(t)$ and $\b C(t)$ is contained in the vertical 
cylinder containing $\b D(q')$. As 
$t \rightarrow 0$, $x_3\res_{C(t)}\rightarrow 0$. Since $\S_1$ is a graph 
in this cylinder, $C(t)$ cannot meet $\S_1$ for the first time at a point 
of $\b C(t)$ (where $C(t)$ is vertical). Also, $C(t)$ cannot touch $\S_1$ at 
an interior point by the maximum principle, nor at infinity. So $C(t)$ 
never touches $\S_1$ and $q_{\infty}$ cannot be in the asymptotic boundary 
of $\S_1$. This proves Theorem~\ref{the9}.
\hfill\qed

\specialnumber{10} \proclaim{Lemma}
\label{catenoid}
Let $C$ be a circle in $\left\{ x_3=0 \right\}$ with center 
$q_{\infty}=(0,0)$. There is a one-parameter family of catenoid cousin 
{\rm (}\/and horosphere\/{\rm )} ends $C(t)$, $0<t\leq 1${\rm ,} satisfying\/{\rm : }
\begin{itemize}
\item[{\rm a)}] each $C(t)$ is a vertical graph over 
$\left\{ 0<x^2+y^2<A^2 \right\}${\rm ,} $A$ the radius of~$C${\rm ,} 
\item[{\rm b)}] $C(t)$ is vertical over $\left\{ x^2+y^2=A^2 \right\}${\rm ,} 
\item[{\rm c)}] $x_3(\b C(1))=A$, $C(1)$ is a horosphere{\rm ,} 
\item[{\rm d)}] $q_{\infty}=\b_{\infty}C(t)$, for each $t${\rm ,} and 
\item[{\rm e)}] $x_3(C(t)) \rightarrow 0$ as $t \rightarrow 0${\rm ,} 
\item[{\rm f)}] $\veh(C(t)). \overrightarrow{e_3} \geq 0${\rm .}
\end{itemize}

\endproclaim

{\it Proof}.
J.M.~Gomes has proved that a family of this nature exists as graphs over 
the exterior domain of $C$ \cite{Go}. To get the $C(t)$ of the lemma, one 
does inversion of this family through a plane $P$ with $\b_{\infty}P=C$, 
followed by a homothety from $q_{\infty}$; cf.\ Figure~14; 
the homothety takes $B$ to $A$. In the appendix we show how these 
surfaces can be obtained.
\hfill\qed
\vglue-12pt
\figin{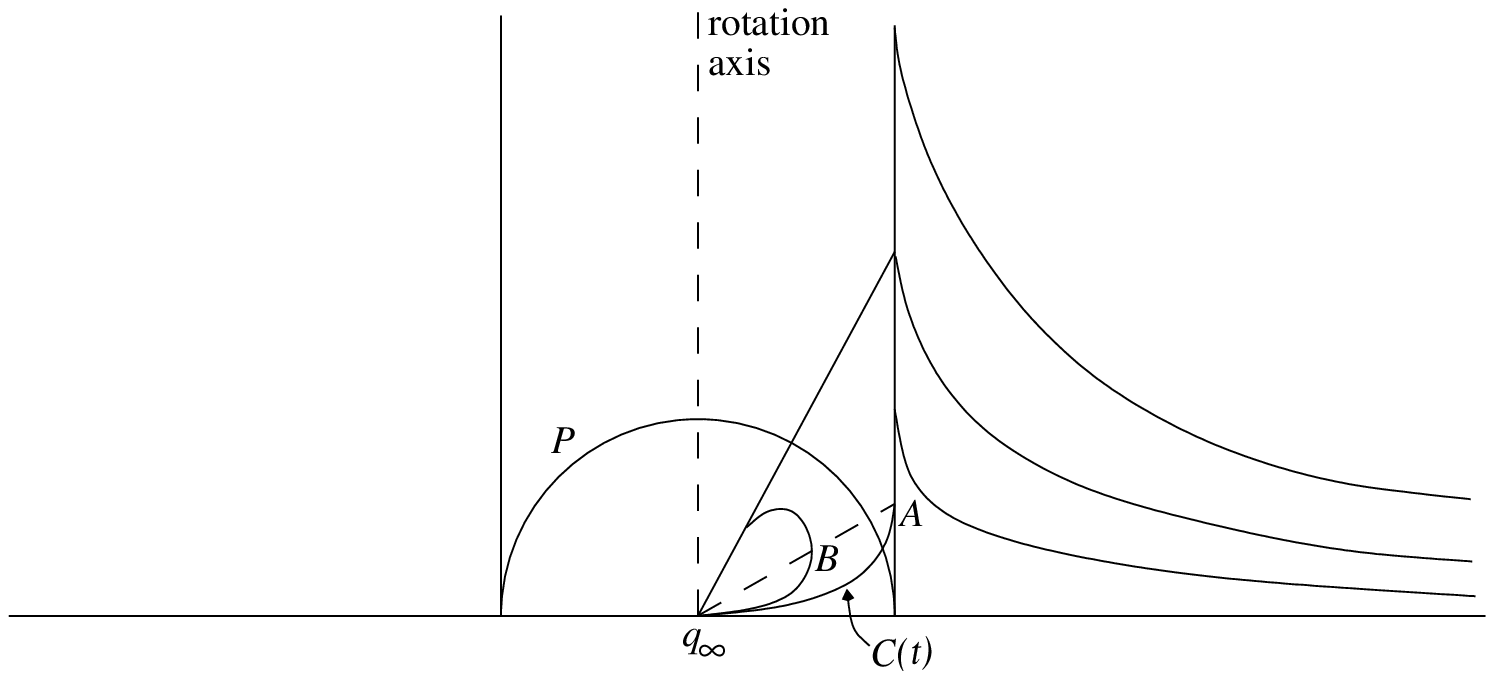}{680}
\centerline{Figure~14}

\specialnumber{10} \proclaim{Theorem}
\label{the10}
Let $E$ be a properly embedded Bryant annular end{\rm .} Then $E$ is not dense 
at infinity{\rm ,} has finite total curvature and is regular{\rm .} Hence {\rm (}\/by 
Corollary~{\rm \ref{standardends})} $E$ is asymptotic to a catenoid cousin end 
or to a horosphere end{\rm .} 
\endproclaim

\vglue-9pt
{\it Proof}.
We remark that the proof of Theorem~\ref{the9} proves Theorem~\ref{the10} 
when $E$ is part of a properly embedded surface $M$ as in 
Theorem~\ref{the9}. Here is the argument in general.

Let $\S$ be a compact embedded surface such that $ \b \S=\b E$ and 
$M=\S \cup E$ is an embedded surface (not necessarily smooth along 
$\b E$). Change the metric of $\h$ in a compact neighborhood of $\S$ so 
that $M$ has mean curvature greater than 1 near $\S$. Now prove 
Theorem~\ref{the8} with $M$ in this new metric. The $\S_1$, $\S_2$ 
one obtains will satisfy all the conditions necessary to do the argument  
of Theorem~\ref{the9}. What matters is the structure of $\S_1$, $\S_2$ 
near infinity. The same argument as in the proof of Theorem~\ref{the9} then 
shows $E$ cannot be dense at infinity. Thus, Corollary~\ref{standardends} 
yields Theorem~\ref{the10}.
\hfill\qed

\vglue-6pt
\specialnumber{11} \proclaim{Theorem}
\label{the11}
Let $M$ be a properly embedded finite topology{\rm ,} Bryant surface{\rm .} If $M$ 
is simply connected {\rm (}\/more generally if $M$ has only one end\/{\rm ),} $M$ is a 
horosphere{\rm .} If $M$ has two ends then $M$ is a catenoid cousin{\rm .} If $M$ 
has three ends{\rm ,} then $M$ is a bigraph over a plane $P${\rm ;} i.e.{\rm ,} $M$ is 
invariant by symmetry in $P$ and each component of $M-P$ is a geodesic 
graph over $P${\rm .}
\endproclaim

\vglue-6pt
\demo{Proof}
When $M$ is simply connected, $\b_{\infty}M$ is one point by 
Theorem~\ref{the10}.  Then M.~do~Carmo and B.~Lawson \cite{docarmo-lawson} 
proved $M$ is a horosphere. When $M$ has two ends, $\b_{\infty}M$ is 
two points and $M$ is invariant by rotations about the geodesic joining 
the two points \cite{L-R}. Thus $M$ is a catenoid cousin. When $M$ has 
three ends, $\b_{\infty}M$ consists of three points so that $\b_{\infty}M$ 
is contained in a circle of $S_{\infty}$. The conclusion is then proved 
in \cite{L-R}.
\enddemo
\vglue-6pt
\specialnumber{12} \proclaim{Theorem}
\label{the12} \hglue-8pt
Let $M$ be a properly embedded Bryant surface{\rm ,} $M$ not a horosphere{\rm .} Then 
each annular end of $M$ is asymptotic to a catenoid cousin~end{\rm .}
\endproclaim

 \vglue-6pt
{\it Proof}.
We know by Theorem~\ref{the10}, that each annular end $E$ is asymptotic 
to a catenoid end or to a horosphere end. We will assume $E$ is asymptotic 
to a horosphere end and obtain a contradiction.

We work in the upper  half-space model of $\HY^3$, $\{ x_3>0 \}$, and 
assume $E$ is asymptotic to a horosphere $x_3=c>0$. In particular the mean 
curvature vector of $E$ points up outside of some compact set of $E$. There 
are no ends of $M$ above $E$ since their mean curvature vector would also 
point up (each such end is asymptotic to a horizontal horosphere or a 
catenoid cousin end whose limiting normal points vertically up) and $M$ 
separates $\HY^3$ into two connected components so that no such end is above $E$.

Then for $\varepsilon>0$, the part $A$ of $M$ above $c+\varepsilon$ is 
compact. At the highest point of $A$  (if $A$ is not empty) the mean 
curvature vector of $M$ points down. But this highest point can be joined 
by an arc in $\HY^3-M$  to a point of $E$ where the mean curvature vector 
points up. Thus $M$ is completely below $x_3=c$.

Let $\varepsilon>0$ and let $C$ be a small circle in the plane 
$x_3=c-\varepsilon$ so that $C$ is above $M$. Just as in the proof of the 
half-space theorem for properly immersed minimal surfaces in $\HY^3$ 
\cite{RR}, one can take a family of catenoid cousin ends $C(\lambda)$, 
$\partial C(1)=C$ with $C(1)$ above $M$, where $C( \lambda )$ converges 
to the plane $x_3=c-\varepsilon$ as $\lambda \rightarrow 0$. Then 
some $C( \lambda )$   touches $M$ at a point $q \in M$ and the maximum 
principle would yield $M$ equals this catenoid cousin. Thus each end of $M$ is 
asymptotic to a catenoid cousin.
\hfill\qed

  \vglue24pt
\centerline{\bf Appendix: The family of graphs of Lemma~\ref{catenoid}}

\bigbreak

Consider the family of vertical catenoids in $\R^3$ whose waist circle is 
of length $\vert\lambda\vert$ and in the $\{x_3=0\}$ plane. Orient by the 
inner pointing normal. The Weierstrass data on the simply connected 
covering space $\C$ are given by $g(z)=e^z$, 
$\omega(z)=\vert\lambda\vert e^{-z}dz$, and the metric is 
$ds=\vert\omega\vert\left( 1+\vert g\vert^2 \right)=
2\vert\lambda\vert\cosh(x)\vert dz\vert$, $z=x+iy$.

The cousins of these catenoids (as $\lambda$ varies) have second 
fundamental form $\widetilde{I\kern -.2em I}=I\kern -.2em I+ds^2$,  and
$I\kern -.2em I$ is  the second fundamental form of the catenoid in $\R^3$. 
The second fundamental form of the catenoid is calculated with respect 
to the inner pointing normal if $\lambda>0$ and the outer normal for 
$\lambda<0$.

One can explicitly find the cousins by solving for $F$ in 
$$F^{-1}dF=\pmatrix{ g&-g^2 \cr 1&-g }\omega.$$
This is done in \cite{Yam-Ume2} and \cite{Rosenberg}, and one obtains in 
the upper half-space model:
\[ \begin{array}{rl}
(x_1+ix_2)(z)=&
\displaystyle{\left[
{\left(\frac{1}{4}-\alpha^2\right)\left(e^x+e^{-x}\right)e^{2\alpha x} 
\over \left(\frac{1}{2}-\alpha\right)^2e^{-x}+
\left(\frac{1}{2}+\alpha\right)^2e^x}\right] e^{2i\alpha y}} \\
\ &\\
x_3(z)=&\displaystyle{2\alpha e^{2\alpha x}\over
\left(\frac{1}{2}-\alpha\right)^2e^{-x}+
\left(\frac{1}{2}+\alpha\right)^2e^x}
\end{array} \]
where $\alpha^2=\frac{1}{4}+\lambda$. This is a surface of revolution for 
$\lambda>-\frac{1}{4}$, embedded for $\lambda>0$ and immersed for 
$-\frac{1}{4}<\lambda<0$.

Let $a=\frac{1}{2}+\alpha$, $b=\frac{1}{2}-\alpha$. The generatrix
$\Gamma$ in the $(x_1,x_3)$ plane of these surfaces of revolution is then 
\eject
\begin{eqnarray*}
\noalign{\vskip-12pt}
x_1(t)&=&
\displaystyle{ab\left(e^t+e^{-t}\right)e^{2\alpha t} 
\over b^2e^{-t}+a^2e^t} ,\\
\noalign{\vskip5pt}
x_3(t)&=&\displaystyle{2\alpha e^{2\alpha t}\over
b^2e^{-t}+a^2e^t}.
\end{eqnarray*}

The points of $\Gamma$ with vertical tangents are the solutions of 
$\frac{\b x_1}{\b t}=0$ and are the solutions of
$$a^2e^{2t}-(2ab+1)+b^2e^{-2t}=0.$$
The discriminent is then $\delta^2=2\left(1-2\alpha^2\right)$, so for 
$0<\alpha<\frac{1}{\sqrt{2}}$ there are two distinct  roots $e^{2\tau}$, 
$e^{2\tau'}$ and $e^{2\tau}e^{2\tau'}=\frac{b^2}{a^2}<1$. We take 
$\tau<\tau'$, so that $e^{2\tau}={2ab+1-\delta\over2a^2}$, and $\tau<0$.

\smallbreak
For $\frac{1}{2}<\alpha<\frac{1}{\sqrt{2}}$, one obtains an embedded 
surface and $\tau=\tau'$ for $\alpha=\frac{1}{\sqrt{2}}$; see Figures~15-a and 15-b. For 
$0<\alpha<\frac{1}{2}$ ($\lambda<0$), one obtains an immersed surface; 
cf.\ Figure~15-c.

\figin{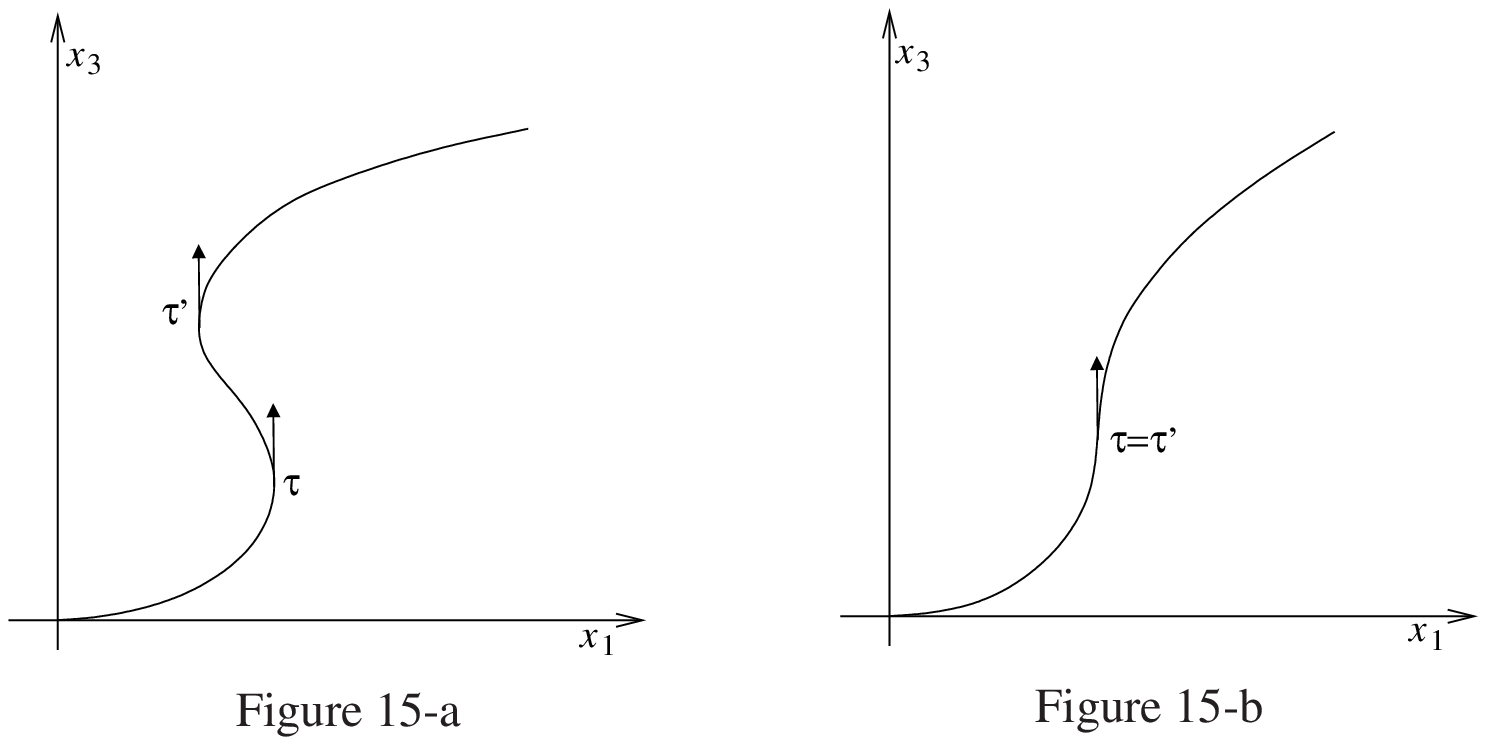}{700}
\figin{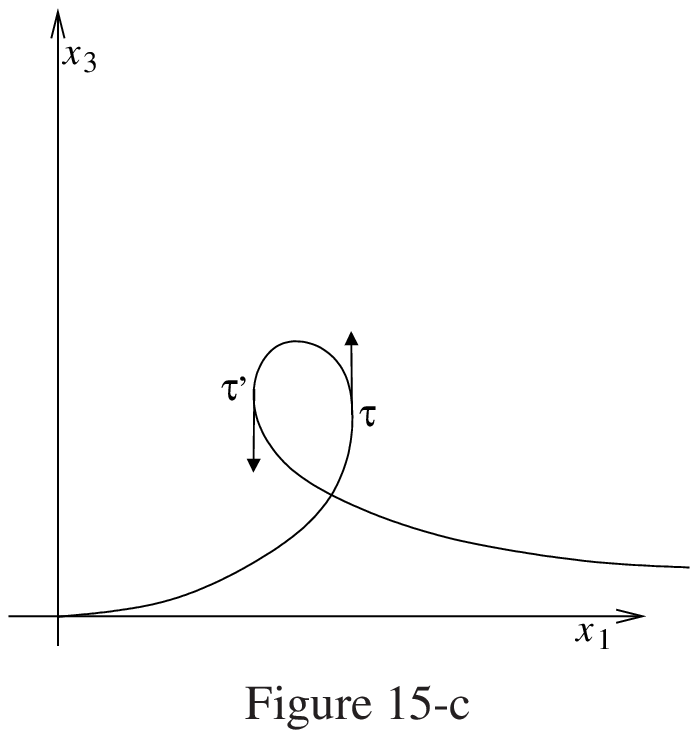}{700}

\vglue-12pt
We are interested in the case $0<\alpha<\frac{1}{2}$. For $-\infty<t<\tau$, 
$\Gamma$ is a graph over an interval $\left(0,x_1(\tau)\right)$. Since we 
want a graph over a fixed interval $\left(0,A\right)$, we renormalize 
by a hyperbolic isometry which is homothety from the origin.

More precisely, let $G_\alpha$ be the graph over $\left(0,A\right)$, 
defined for $-\infty<t<\tau$. We have on $G_\alpha$:
\begin{eqnarray*}
x_1(t)&=&
\displaystyle{{A\over x_1(\tau)}{ab\left(e^t+e^{-t}\right)e^{2\alpha t} 
\over b^2e^{-t}+a^2e^t}} ,\\
\noalign{\vskip4pt}
x_3(t)&=&\displaystyle{{A\over x_1(\tau)}{2\alpha e^{2\alpha t}\over
b^2e^{-t}+a^2e^t}}\ .
\end{eqnarray*}
Hence 
$${x_3(t)\over x_1(t)}={\alpha \over ab\cosh(t)}\leq
{\alpha \over ab\cosh(\tau)},$$ since $t<\tau<0$. It is easy to 
see that $\displaystyle\lim_{\alpha\rightarrow0}
\left(\alpha\over ab\cosh(\tau)\right)=0$; hence the graphs limit to 
$\left(0,A\right)\times\{0\}$ as $\alpha\rightarrow0$, as desired.

\end{document}